\documentclass[11pt,american]{article}
\usepackage[T1]{fontenc}
\usepackage[utf8]{inputenc}
\usepackage{geometry}
\usepackage{color}
\usepackage{babel}
\usepackage{float}
\usepackage{amsmath}
\usepackage{amsthm}
\usepackage{amssymb}
\usepackage{graphicx}
\usepackage{esint}
\usepackage[unicode=true,pdfusetitle,
 bookmarks=true,bookmarksnumbered=false,bookmarksopen=false,
 breaklinks=false,pdfborder={0 0 0},backref=false,colorlinks=true]
 {hyperref}

\makeatletter

\floatstyle{ruled}
\newfloat{algorithm}{tbp}{loa}
\providecommand{\algorithmname}{Algorithm}
\floatname{algorithm}{\protect\algorithmname}

\numberwithin{equation}{section}
\theoremstyle{plain}
\newtheorem{thm}{\protect\theoremname}[section]
  \theoremstyle{definition}
  \newtheorem{defn}[thm]{\protect\definitionname}
  \theoremstyle{remark}
  \newtheorem{rem}[thm]{\protect\remarkname}
  \theoremstyle{definition}
  \newtheorem{problem}[thm]{\protect\problemname}

\usepackage{algorithm,algpseudocode}
\usepackage{url}

\@ifundefined{showcaptionsetup}{}{%
 \PassOptionsToPackage{caption=false}{subfig}}
\usepackage{subfig}
\makeatother

  \providecommand{\definitionname}{Definition}
  \providecommand{\problemname}{Problem}
  \providecommand{\remarkname}{Remark}
\providecommand{\theoremname}{Theorem}

\begin{document}

\title{A modified algebraic reconstruction technique taking refraction into
account with an application in terahertz tomography}

\author{Jens Tepe\thanks{Department of Mathematics, Saarland University, PO Box 15 11 50, 66041
Saarbr\"ucken, Germany (jens.tepe@num.uni-sb.de)},\; Thomas Schuster\thanks{Department of Mathematics, Saarland University, PO Box 15 11 50, 66041
Saarbr\"ucken, Germany (thomas.schuster@num.uni-sb.de), correspondent
author.}\; and Benjamin Littau\thanks{German Plastics Center, 97076 W\"urzburg, Germany (b.littau@skz.de)}}

\date{~\vspace{-1cm}
}

\maketitle


\section*{Abstract}

{\footnotesize{}Terahertz (THz) tomography is a rather novel technique
for nondestructive testing that is particularly suited for the testing
of plastics and ceramics. Previous publications showed a large variety
of conventional algorithms adapted from computed tomography or ultrasound
tomography which were directly applied to THz tomography. Conventional
algorithms neglect the specific nature of THz radiation, i.~e.~refraction
at interfaces, reflection losses and the beam profile (Gaussian beam),
which results in poor reconstructions.  The aim is the efficient
reconstruction of the complex refractive index, since it indicates inhomogeneities
in the material. A hybrid algorithm has been developed based on the
algebraic reconstruction technique (ART). ART is adapted by including
refraction (Snell's law) and reflection losses (Fresnel equations). Our method
uses a priori information about the interface and layer geometry
of the sample. This results in the ``Modified ART for THz tomography'',
which reconstructs simultaneously the complex refractive index from
transmission coefficient and travel time measurements.}{\footnotesize \par}

\section{Introduction}

Terahertz (THz) radiation has become increasingly popular for nondestructive
testing of plastics and ceramics in recent years. With a frequency
range between $100\,\mathrm{GHz}$ to $10\,\mathrm{THz}$, this radiation
is situated between microwaves and infrared light within the electromagnetic
spectrum. Thus, in terms of the wavelength, THz radiation shows elements
of ray character (X-ray CT) as well as of wave character (ultrasound
tomography). As a result,  the wave character cannot be completely
neglected. While plastics and ceramics are almost transparent in
this frequency range, THz radiation cannot easily penetrate water,
metals and conductive materials. Because of its relatively long wave
length compared to X-rays, THz radiation has a limited lateral
resolution. However, by using modulated THz radiation, information
about both phase and amplitude are available as measured data in contrast
to conventional tomographic methods \cite{Krumbholz2009}. 

The inverse problem of THz tomography consists of determining dielectric
properties (i.~e.~refractive index $n$ and absorption coefficient
$\alpha$) of the material from phase and amplitude information of
the THz radiation, which admits finding material characteristics such
as filler content or the detection of cracks and air cavities. For
an overview on the relatively new field of tomographic methods using
THz radiation, see Wang et~al.~\cite{Wang2004}, Guillet et~al.~\cite{Guillet2014}
and Ferguson \cite{Ferguson2004,Ferguson2002}. Published results
can be divided into transmission or reflection tomographic methods
based on the available measured data. The adaption of computed tomography
methods \cite{Kak1979,Natterer2001} to the THz range is widely distributed
\cite{Brahm2010,Ewers2009,Ferguson2002,Jewariya2013,Nguyen2006,Wang2004}.
Here a very simplified propagation along straight lines is used with
the Radon transform as mathematical model. However, physical effects
such as refraction, diffraction, edge effects and the Gaussian beam
profile have an influence on THz measurements \cite{Brahm2014}. Some
publications consider these effects by adaptations. In \cite{Recur2012},
Recur et~al.~presented a filtered backprojection and adapted iterative
methods (SART, OSEM) for 3D THz CT, which considers the intensity
distribution of THz rays (Gaussian beam). Mukherjee et~al.~only
regard cylindric objects, but take Fresnel reflection losses and beam
steering losses into account, see \cite{Mukherjee2013}. Furthermore,
Ferguson et~al.~\cite{Ferguson2002} and Wang et~al.~\cite{Wang2004}
reconstruct the refractive index using measurements of the diffracted
THz field based on the Helmholtz equation (so-called T-ray Diffraction
Tomography, analogous to ultrasound tomography, see \cite{Colton2000,Colton2013}). 

Due to the long wavelength of THz radiation compared to X-rays (for
a frequency of $f=0.1\,\mathrm{THz}$, the wavelength is approx. $\lambda=3\,\mathrm{mm}$),
wave phenomena need to be taken into account in the modeling of the
measurement process. On the other hand, a more precise model leads
to a computationally more intensive solution. Our task is to find
a balance between both. 

In the following we present a novel, hybrid algorithm for THz tomography.
The algorithm is developed for a tomograph designed for medium-sized
companies, which demands a short runtime. Because of that we search
for an efficient solution of the inverse problem by considering only
dominating effects. 

\begin{figure}
\noindent \hfill{}\subfloat[3D representation of the Gaussian beam \label{fig:3D-presentation-Gaussian-Beam} ]{

\includegraphics[width=0.45\textwidth]{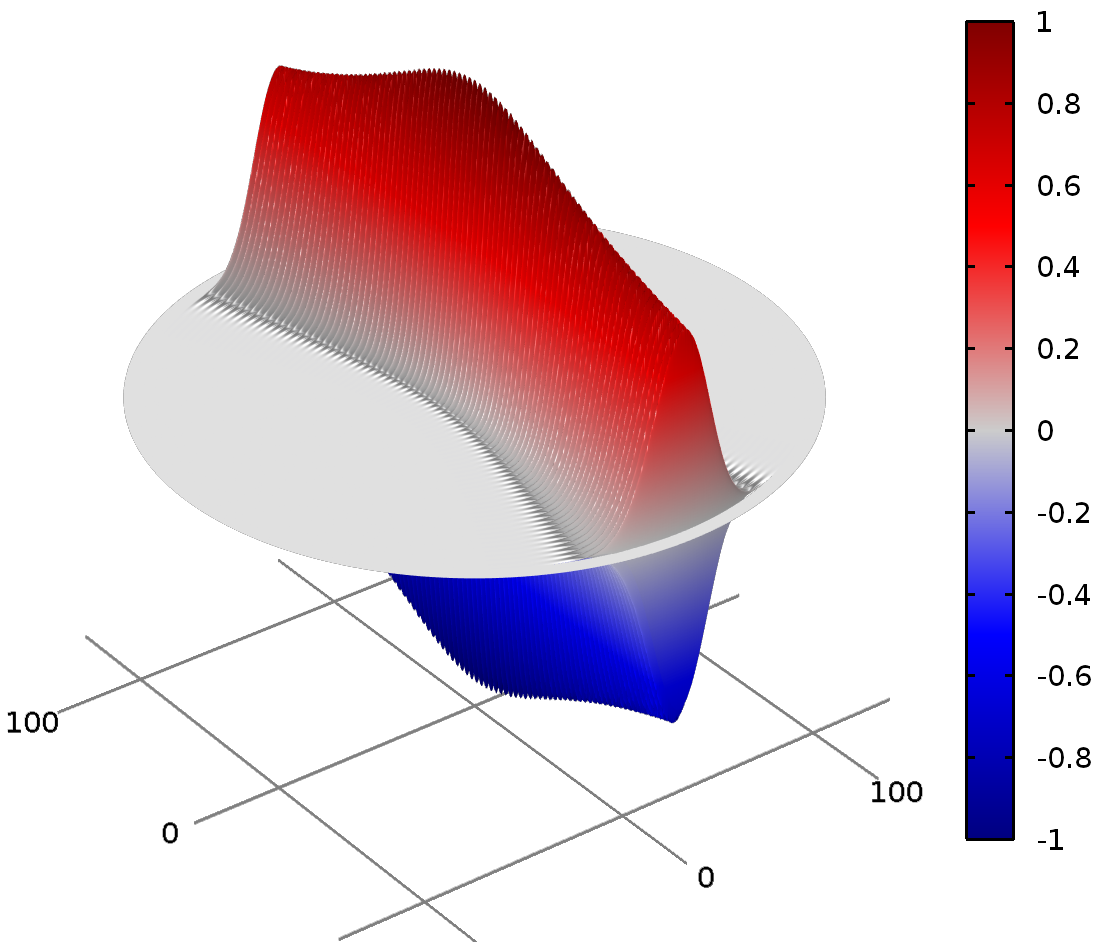}}\hfill{}\subfloat[Simulation of the electric field of an object consisting of polyethylene
($\tilde{n}=1.51+0.023i$) \label{fig:Helmholtz-broken-ray}]{

\includegraphics[width=0.5\textwidth]{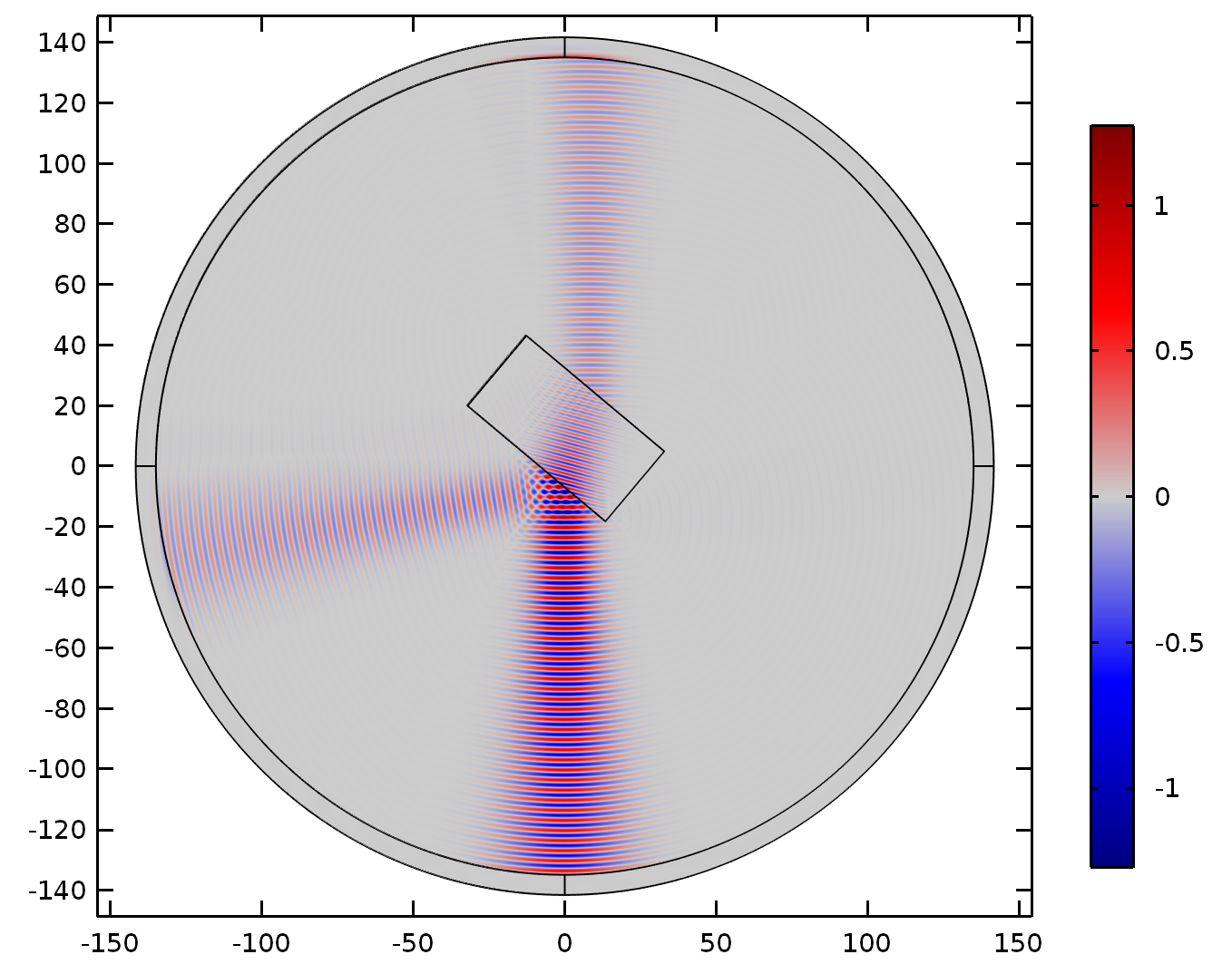}}\hfill{}\caption{THz radiation for $f=90\,\mathrm{GHz}$ propagating from bottom to
top as a simulation based on the Helmholtz equation \label{fig:THz-radiation-simulation}}
\end{figure}

The radiation from a THz source behaves ideally as a Gaussian beam
and can be modeled by the wave equation or, if we are only interested
in the intensity, by the Helmholtz equation. The transverse intensity
profile is well approximated by Gaussian functions (see figure \ref{fig:3D-presentation-Gaussian-Beam}).
However, the solution of the Helmholtz equation is very time-consuming.
The wave length is significantly smaller than in the case of acoustic
waves, forcing a fine discretization of the domain of interest. This
prevents an efficient solution with ultrasound tomography algorithms.
Hence, our goal is a new algorithm not based on the solution of wave
or Helmholtz equation, but still taking into account the main physical
effects. These are primarily transmission, refraction at interfaces,
reflection losses and the Gaussian beam profile (see figure \ref{fig:Helmholtz-broken-ray}).
 In the following we will neglect the Gaussian beam profile for the
sake of simplicity.

The THz tomography system used for this work provides intensity data,
analogous to the conventional transmission tomography, as well as
(delayed) travel times or path differences, as in ultrasound technology.
The reconstruction algorithm described here uses travel time and transmission
data simultaneously to reconstruct the material parameters refractive
index and absorption coefficient. Therefore, we adapt the algebraic
reconstruction technique (ART, see \cite{Gordon1970}) to our physical
conditions by taking into account refraction at the interfaces of
the object according to Snell's law and including the reflection losses
given by the Fresnel equations. For this, we need some a priori information
on the object's interfaces. 

This approach is comparable to those in ultrasound CT and ultrasound
vector field tomography, where reconstructions from travel time measurements
are improved by considering ray paths which take deflection of the
ray into account. In ultrasound CT, the ray paths are computed by
solving the eikonal equation (ray tracing) \cite{Andersen1982,Schomberg1978,Smith1980}.
This approach unfortunately requires a continuously differentiable
refractive index and thus a smoothing of any discontinuities of the
refractive index at interfaces. Furthermore, it is necessary to solve
a nonlinear partial differential equation of first order. Ray tracing
algorithms were successfully combined with algebraic reconstruction
techniques (ART, SART) \cite{Andersen1987,Andersen1991,Denis1995}.
These algorithms work well, if discontinuities are smaller than $10-20\%$
refractive index deviation \cite{Andersen1990,Denis1995}. In THz
tomography, these deviations at interfaces are even higher compared
to ultrasound tomography (about $50\%$). Therefore, refraction at
interfaces can be assumed as a dominating effect in THz tomography. 

In ultrasound vector field tomography, signals travel along geodesic
curves of a Riemannian metric, which is associated with the refractive
index according to Fermat's principle \cite{Schroeder2015,Pfitzenreiter2011}.
The inverse problem consists of determining the refractive index from
integrals along geodesics curves from ultrasound time-of-flight (TOF)
measurements. In both, THz tomography and ultrasound vector field
tomography, the ray path resp.~geodesic curve depends on the unknown
refractive index.

 The algorithm presented here allows discontinuities in the refractive
index and includes the change of direction of the rays at interfaces
by a direct application of Snell's law. 

The paper is organized as follows. In the next section we briefly
recall the mathematical model for the 2D THz transmission and travel
time tomography. Subsequently, the notation of the a priori given
interfaces and a brief description of the physical model, especially
the Fresnel equations and Snell's law, will be discussed. In Section
\ref{sub:ray-path-depending-on-refraction} we derive a description
of ray paths considering refraction. Section \ref{sec:algorithm-development}
deals with the development of the new algorithm, beginning with a
brief revision of the conventional ART. Based on ART, we introduce
the modified ART for THz tomography and its implementation. Eventually,
the validation of the modified ART and its implementation are to be
found in section \ref{sec:Numerical-experiments} using synthetic
as well as real measured data, comparing our algorithm with the filtered
backprojection (FBP) and conventional ART. We conclude the paper with
a summary of the main results and an outlook for future work and potential
improvements.

\section{The mathematical model of 2D THz tomography\label{sec:Mathematical-model}}

We consider 2D THz transmission and travel time tomography with standard
parallel scanning geometry. The main task is the reconstruction of
an object characterized by the complex refractive index $\tilde{n}:\Omega\rightarrow\mathbb{C}$
with
\begin{equation}
\tilde{n}(x)=n(x)+i\kappa(x),\qquad x\in\Omega,
\end{equation}
where $\Omega=\left\{ x\in\mathbb{R}^{2}:\left\Vert x\right\Vert <R\right\} $
with $R>0$ the radius of the considered domain. The real part $n$
describes the optical path length of the material, whereas the imaginary
part $\kappa$ models the material-dependent absorption. The transmission
coefficient $\tau=I/I_{0}$ and the path difference $d=c_{0}(T-T_{0})$
are available as measured data, with $I_{0}$ the initial intensity,
$I$ the intensity at the detector, $T$ the travel time with and
$T_{0}$ without an object in the ray path and $c_{0}$ the speed
of light. The real part of the complex refractive index is the quotient
of the speed of light $c_{0}$ and the propagation speed  $c_{M}$
inside the medium, 

\begin{equation}
n(x)=\frac{c_{0}}{c_{M}(x)},\qquad x\in\Omega,\,c_{M}>0.
\end{equation}
This means the (delayed) travel time is physically related to real
part of $\tilde{n}$. The imaginary part 

\begin{equation}
\kappa(x)=\alpha(x)\frac{c_{0}}{4\pi f},\qquad x\in\Omega,\,f>0,
\end{equation}
describes the absorption of the wave by the absorption coefficient
$\alpha$ in $\mathrm{cm^{-1}}$ and is connected to the transmission
coefficient $\tau$. Furthermore, we assume that the complex refractive
index vanishes outside the reconstruction region $\Omega$, i.~e.~$n(x)-1=0$
and $\kappa(x)=0$, if $x\notin\Omega$. Let $L$ be the propagation
path of the THz radiation. The ratio of the transmitted intensity
$I$ and the initial intensity $I_{0}$ is connected to the absorption
coefficient $\alpha$ via 

\begin{equation}
g_{\mathrm{abs}}(L)=\intop_{L}\alpha(x)\,\mathrm{d}x=\ln\left(\frac{I_{0}}{I}\right)=\ln\left(\frac{1}{\tau}\right)\label{eq:integral-eq-CT}
\end{equation}
according to Lambert-Beer's law. The path difference is given by

\begin{equation}
g_{\mathrm{ref}}(L)=\intop_{L}\left(n(x)-1\right)\,\mathrm{d}x=d,
\end{equation}
which follows from subtraction of the path in air with a refractive
index of $n_{0}=1$, see \cite{Natterer1986}. 

Due to refraction at interfaces, THz rays do not necessarily propagate
along straight lines. The rays are rather refracted at interfaces
according to Snell's law. In case of a piecewise constant refractive
index, the propagation path is continuous and consists of piecewise
straight lines. We name the segments of a THz ray between two interfaces
as ``partial ray'' $L_{l}\subseteq L$ with $l=1,...,\hat{K}+1$ and
$\hat{K}$ the number of intersected interfaces. The first partial
ray $L_{1}$ is parametrized as in X-ray CT until reaching the first
interface, see \cite[p. 41-42]{Natterer2001}. Let 
\begin{equation}
\Gamma:\mathbb{R}\rightarrow\mathbb{R}^{2},\qquad\Gamma(t)=s\omega(\varphi)+t\omega^{\bot}(\varphi)
\end{equation}
be the parametrization of $L_{1}$ with distance
$s\in\mathbb{R}$ from the origin, the angle $\varphi\in\left[0,\,2\pi\right)$, $\omega(\varphi)=\left(\cos\varphi,\sin\varphi\right)^{T}$
and $\omega^{\bot}=(-\sin\varphi,\,\cos\varphi)^{T}$. Using $\Gamma$,
the ray is represented by

\begin{equation}
L_{1}\left(\varphi,s\right)=\left\{ \Gamma(t):\,t\in\left[-\infty,\,t_{1}^{\mathrm{if}}\left(\varphi,s\right)\right]\right\} \label{eq:parametrization-L}
\end{equation}
with $t_{1}^{\mathrm{if}}\left(\varphi,s\right)$ the time of the
incidence on the first interface (so-called ``initial ray''). In
X-ray CT, this parameter is usually set to $t_{1}^{\mathrm{if}}=\infty$.
The full ray path depending on refraction is described in section
\ref{sub:ray-path-depending-on-refraction}. Therefore we first need
to take a closer look at a suitable notation of the interfaces and
afterwards we briefly describe the relevant physical laws.

\subsection{Notation of interfaces}

\begin{figure}[H]
\begin{centering}
\includegraphics[width=0.5\textwidth]{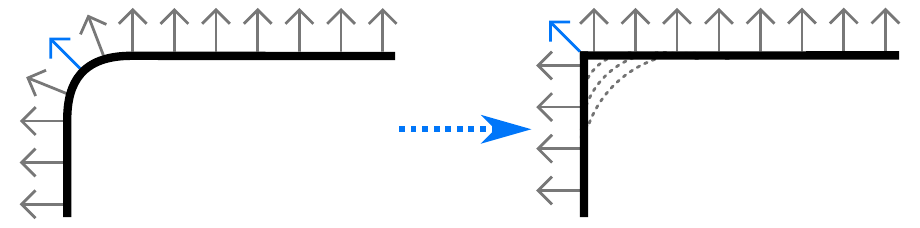}
\par\end{centering}

\caption{Dealing with corners: the calculated normal for the discontinuity
corresponds to the normal of the smoothed interface at this point\label{fig:Dealing-with-corners}}
\end{figure}

As mentioned before, we would like to use some a priori information.
We postulate the knowledge of the object's interfaces, i.~e.~our
a priori information contains positions and the normals of the interfaces
of the object. This is quite a strong assumption, but in many cases
the (approximate) location of some interfaces is actually known, for
example from CAD drawings of the examined object or preliminary testings
with an X-ray CT. However, the interfaces of potential defects are
unknown. Nevertheless, the knowledge of some interfaces already leads
to better reconstructions. 
\begin{defn}
An \emph{interface} is a regular continuously differentiable and injective
immersion $\text{\ensuremath{\Xi}}:\left[a,b\right]\rightarrow\Omega\subset\mathbb{R}^{2}$,
fulfilling $\Xi^{'}(\sigma)=\frac{d\Xi}{d\sigma}\neq0,\,\forall\sigma\in\left[a,b\right]$.

\end{defn}
The curve has a tangent and a normal $\mathfrak{n}$ at every point,
because the first derivative $\Xi^{'}$ is continuous and does not
vanish. The tangent vector of $\Xi$ at $\sigma_{0}\in\left[a,\,b\right]$
is given by
\begin{equation}
\Xi^{'}(\sigma_{0})=\left.\frac{d\Xi}{d\sigma}\right|_{\sigma=\sigma_{0}}=\left(\xi_{1}^{'}(\sigma_{0}),\,\xi_{2}^{'}(\sigma_{0})\right)^{T},
\end{equation}
whereby we define the normal vector 

\begin{equation}
\mathfrak{n}(\sigma_{0})=\left(-\xi_{2}^{'}(\sigma_{0}),\,\xi_{1}^{'}(\sigma_{0})\right)^{T},
\end{equation}
which we will need for Snell's law. Let $K\in\mathbb{N}$ be the number
of known interfaces $\Xi_{k},\,k=1,...,\,K$ of the object. Interfaces
may touch or intersect each other as we often have to deal with non-smooth
interfaces, e. g. rectangular objects. However, $\Xi^{'}(\sigma)$
and hence the normal does not exist at corners. This case is handled
by local smoothing. Let us assume, that
\begin{equation}
P:=\bigcap_{\iota=1}^{\tilde{K}}\Xi_{k_{\iota}}
\end{equation}
is the shared intersection of $\tilde{K}\in\left[1,\,K\right]$ different
interfaces $\Xi_{k_{\iota}}$ with $k_{\iota}\in\left[1,\,K\right],\iota=1,...,\,\tilde{K}$.
Then $\sigma_{k_{\iota}}:=\Xi_{k_{\iota}}^{-1}\left(P\right)\in\left[a_{k_{\iota}},b_{k_{\iota}}\right]$
parametrizes the point of intersection on each interface and $\mathfrak{n}_{k_{\iota}}(\sigma_{k_{\iota}})$
are the corresponding normals for $\iota=1,...,\tilde{K}$. We now
compute the normal in $P$ by averaging
\begin{equation}
\mathfrak{n}=\frac{\sum_{\iota=1}^{\tilde{K}}\mathfrak{n}_{k_{\iota}}(\sigma_{k_{\iota}})}{\left\Vert \sum_{\iota=1}^{\tilde{K}}\mathfrak{n}_{k_{\iota}}(\sigma_{k_{\iota}})\right\Vert },\label{eq:sum-normals}
\end{equation}
which can be interpreted as a mean orientation of the interfaces containing
$P$. This approach sets a normal for the intersection, which occurs
by a local smoothing of the corner $P$ (see figure~\ref{fig:Dealing-with-corners}).
For $\tilde{K}=1$ the normal vector is just normalized.

\subsection{Refraction and reflection losses}

We follow the idea of Mukherjee et. al, see \cite{Mukherjee2013},
who consider refraction and reflection losses for cylindrical objects
with known refractive index for the reconstruction of the absorption
coefficient. We generalize this approach for general objects with
given interfaces, unknown $n$ and reconstruct both refractive index
and absorption coefficient. The propagation of electromagnetic waves
at an interface between two different media is influenced by refraction,
reflection, diffraction and scattering. As stated before, we want
to neglect diffraction and scattering. Let us consider the transition
of a ray from a medium with refractive index $n_{1}$ into a medium
with refractive index $n_{2}$. The change in the direction of propagation
of the transmitted part of the THz ray can be calculated according
to Snell's law
\begin{equation}
n_{1}\sin\gamma_{1}=n_{2}\sin\gamma_{2},\label{eq:Snell-1}
\end{equation}
with the angle of incidence $\gamma_{1}$ and the angle of refraction
$\gamma_{2}$ measured from the normal $\mathfrak{n}$ of the interface,
see \cite[p. 38 ff.]{Born1970}. When rays travel from a medium with
a higher refractive index $n_{1}$ to one with a lower refractive
index $n_{2}$, total reflection is possible for a sufficiently large
$\gamma_{1}$. In this case, we have 
\begin{equation}
\frac{n_{1}}{n_{2}}\sin\gamma_{1}>1,\qquad\mathrm{for\,}\,n_{1}>n_{2}
\end{equation}
and the radiation is fully reflected.
\begin{figure}
\begin{centering}
\includegraphics[width=0.5\textwidth]{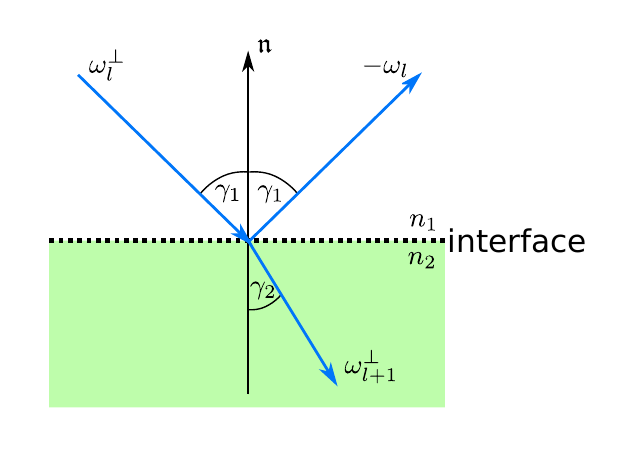}
\par\end{centering}

\caption{Snell's law with normal $\mathfrak{n}$ of the interface and direction
vectors $\omega_{l}^{\perp}$, $\omega_{l+1}^{\perp}$ and $-\omega_{l}$
of the incident, refracted and reflected rays, $n_{1}<n_{2}$\label{fig:draft-Snells-law}}
\end{figure}
The fraction of the incident radiation that is reflected at the interface
can be described by the Fresnel equations, see \cite[p. 38 ff.]{Born1970}.
Reflection losses depend on the polarization of the incident wave,
the angle of incident $\gamma_{1}$ and refraction $\gamma_{2}$ and
the refractive indices. The reflectance of a ray that is polarized
perpendicularly to the plane of incidence is given by

\begin{equation}
\rho=\left|\frac{E_{r}}{E_{0}}\right|^{2}=\left|\frac{n_{1}\cos\gamma_{1}-n_{2}\cos\gamma_{2}}{n_{1}\cos\gamma_{1}+n_{2}\cos\gamma_{2}}\right|^{2},\label{eq:Fresnel}
\end{equation}
where $E_{0}$ resp.~$E_{r}$ is the electric field amplitude of
the incident resp.~reflected partial ray. Analogous considerations
exist also for parallel polarization. All other polarizations can
be represented by a superposition of the linearly perpendicular and
linearly parallel case.

\subsection{Ray propagation considering refraction\label{sub:ray-path-depending-on-refraction}}

\begin{figure}
\begin{centering}
\includegraphics[width=0.5\textwidth]{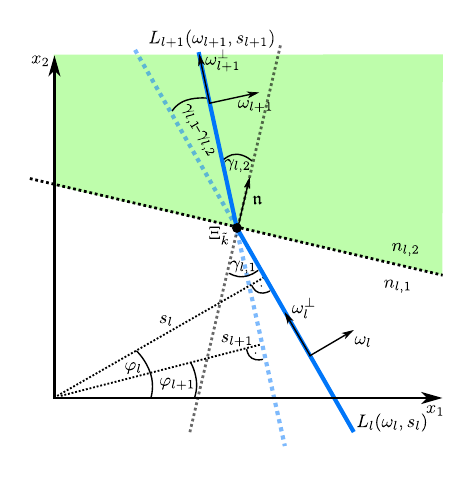}
\par\end{centering}

\caption{Snell's law and parametrization of a refracted THz ray \label{fig:Snell-parametrization}}
\end{figure}

The path of a THz ray $L$ up to reaching the first interface can
be described by equation (\ref{eq:parametrization-L}). The ray is
refracted as it transitions into the next medium and the ray direction
has to be recalculated according to Snell's law. We first have to
determine the next intersection with an interface for each (refracted)
partial ray $L_{l}$. The change in the direction of the ray is calculated
in this intersection. Let the $l$-th partial ray $L_{l}$ be given
by an angle $\varphi_{l}$ and a distance $s_{l}$, i~.e.~$L_{l}\left(\varphi_{l},s_{l}\right)=\left\{ \Gamma_{l}(t)=s_{l}\omega_{l}+t\omega_{l}^{\perp}:\,t\in\left[t_{l-1}^{\mathrm{if}},\,t_{l}^{\mathrm{if}}\right]\right\} $.
Now we calculate the intersections of $L_{l}$ with the interfaces
and determine the parameter 

\begin{equation}
t_{l}^{\mathrm{if}}\left(\varphi_{l},s_{l}\right)=\underset{t}{\min}\left\{ \left(t,\sigma\right)\in\left(t_{l-1}^{\mathrm{if}},\infty\right]\times\left[a_{k},b_{k}\right]:\,\Gamma_{l}(t)=\Xi_{k}(\sigma),\,k=1,...,K\right\} ,\label{eq:t-of-the-next-interface}
\end{equation}
which defines the nearest intersection with an interface in the direction
of propagation with $l\geq1$ and $t_{0}^{\mathrm{if}}=-\infty$.
Let $\Xi_{\tilde{k}}(\sigma),\,\tilde{k}\in\left[1,\,K\right],$ be
the nearest interface and $\tilde{\sigma}$ the parameter corresponding
to $t_{l}^{\mathrm{if}}$ in (\ref{eq:t-of-the-next-interface}) of
the intersection $\Xi_{\tilde{k}}(\tilde{\sigma})=s_{l}\omega_{l}+t_{l}^{\mathrm{if}}\omega_{l}^{\perp}$.
We can determine the refractive indices for both sides of the interface,
which we need for Snell's law, by

\begin{equation}
\begin{array}{lllc}
n_{l,1}:=n\left(\Xi_{\tilde{k}}(\tilde{\sigma})-\varepsilon\mathfrak{n}\right) & \mathrm{\,and\,} & n_{l,2}:=n\left(\Xi_{\tilde{k}}(\tilde{\sigma})+\varepsilon\mathfrak{n}\right), & \mathrm{for\,}\left\langle \omega_{l}^{\perp},\mathfrak{n}\right\rangle \leq0\\
n_{l,1}:=n\left(\Xi_{\tilde{k}}(\tilde{\sigma})+\varepsilon\mathfrak{n}\right) & \mathrm{\,and\,} & n_{l,2}:=n\left(\Xi_{\tilde{k}}(\tilde{\sigma})-\varepsilon\mathfrak{n}\right), & \mathrm{for\,\left\langle \omega_{l}^{\perp},\mathfrak{n}\right\rangle }>0
\end{array}\label{eq:n-for-both-sides-of-the-if}
\end{equation}
with the normal $\mathfrak{n}:=\mathfrak{n}(\tilde{\sigma})$ as in
(\ref{eq:sum-normals}) and an appropriately small $\varepsilon>0$.
The angle of incidence $\gamma_{l,1}$ between the direction vector
of the incident ray $\omega_{l}^{\perp}$ and normal $\mathfrak{n}$
(see figure \ref{fig:draft-Snells-law}) can be calculated by the
scalar product in $\mathbb{R}^{2}$

\begin{equation}
\gamma_{l,1}=\arccos\left|\left\langle \omega_{l}^{\perp},\,\mathfrak{n}\right\rangle \right|
\end{equation}
and the angle of the refracted ray can be obtained by
\begin{equation}
\gamma_{l,2}=\left\{ \begin{array}{lcl}
{\displaystyle \pi-\gamma_{l,1}}, &  & \mathrm{for\,\,}n_{l,1}>n_{l,2}\wedge\gamma_{l,1}\geq\arcsin\frac{n_{l,2}}{n_{l,1}}\mathrm{\,\,(total\,reflection)},\\
{\displaystyle \arcsin\frac{n_{l,1}\sin\gamma_{l,1}}{n_{l,2}},} &  & \mathrm{otherwise}.
\end{array}\right.
\end{equation}
The next direction vector $\omega_{l+1}$ can be easily derived by
a rotation of $\omega_{l}:=\omega\left(\varphi_{l}\right)$ around
an angle $\theta=\pm\left(\gamma_{l,1}-\gamma_{l,2}\right)$ with
the rotation matrix 

\begin{equation}
R_{\theta}:=\left(\begin{array}{cc}
\cos\theta & -\sin\theta\\
\sin\theta & \cos\theta
\end{array}\right)
\end{equation}
see figure \ref{fig:Snell-parametrization}.  The sign of $\theta$
depends on the orientation of the normal $\mathfrak{n}$ to the incident
ray, see (\ref{eq:omega-l-plus-1-cases}). After applying angle sum
identities, the new angle $\varphi_{l+1}$ and the unit vector $\omega_{l+1}$
are computed by

\begin{equation}
\omega_{l+1}\left(\varphi_{l+1}\right):=\left\{ \begin{array}{ll}
\begin{array}[t]{c}
\overset{=:\varphi_{l+1}}{\omega(\overbrace{\begin{array}{llc}
\gamma_{l,1}-\gamma_{l,2} & + & \varphi_{l}\end{array}}}),\end{array} & \mathrm{for\,}\left\langle \omega_{l}^{\perp},\,\mathfrak{n}\right\rangle <0\wedge\left\langle \omega_{l}^{\perp},\,\mathfrak{n}^{\perp}\right\rangle \leq0\\
\begin{array}{c}
\omega(\begin{array}{llc}
\gamma_{l,1}-\gamma_{l,2} & + & \varphi_{l}\end{array}),\end{array} & \mathrm{for\,}\left\langle \omega_{l}^{\perp},\,\mathfrak{n}\right\rangle >0\wedge\left\langle \omega_{l}^{\perp},\,\mathfrak{n}^{\perp}\right\rangle \geq0\\
\begin{array}{c}
\omega(\begin{array}{llc}
\gamma_{l,2}-\gamma_{l,1} & + & \varphi_{l}\end{array}),\end{array} & \mathrm{for\,}\left\langle \omega_{l}^{\perp},\,\mathfrak{n}\right\rangle \leq0\wedge\left\langle \omega_{l}^{\perp},\,\mathfrak{n}^{\perp}\right\rangle >0\\
\begin{array}[b]{c}
\omega(\begin{array}[t]{llc}
\underset{=:\theta}{\underbrace{\gamma_{l,2}-\gamma_{l,1}}} & + & \varphi_{l}\end{array}),\end{array} & \mathrm{for\,}\left\langle \omega_{l}^{\perp},\,\mathfrak{n}\right\rangle \geq0\wedge\left\langle \omega_{l}^{\perp},\,\mathfrak{n}^{\perp}\right\rangle <0
\end{array}\right..\label{eq:omega-l-plus-1-cases}
\end{equation}
The next distance from the origin $s_{l+1}$ of partial ray $L_{l+1}$
follows directly from the intersection of the lines $\left(0+u\omega_{l+1}\right)\cap\left(\Xi_{\tilde{k}}(\tilde{\sigma})+t\omega_{l+1}^{\perp}\right)$
with parameters $t$ and $u$. We obtain with $\Xi_{\tilde{k}}(\tilde{\sigma}):=\left(\xi_{1},\xi_{2}\right)^{T}$ 

\begin{equation}
s_{l+1}=\left|\left(0+u\omega_{l+1}\right)\cap\left(\Xi_{\tilde{k}}(\tilde{\sigma})+t\omega_{l+1}^{\perp}\right)\right|=\left|\cos\varphi_{l+1}\xi_{1}+\sin\varphi_{l+1}\xi_{2}\right|.\label{eq:s-l-plus-1-calculation}
\end{equation}
Consequently, the complete refracted ray path is given by the union
\begin{equation}
L\left(\varphi,s\right):=\bigcup_{l=1}^{\hat{K}+1}\left\{ \Gamma_{l}(t):\,t\in\left[t_{l-1}^{\mathrm{if}},\,t_{l}^{\mathrm{if}}\right]\right\} 
\end{equation}
of all partial rays $L_{l}\left(\varphi_{l},s_{l}\right)$ defined
by (\ref{eq:omega-l-plus-1-cases}) and (\ref{eq:s-l-plus-1-calculation})
with an initial direction $\varphi=\varphi_{1}$ and distance $s=s_{1}$,
$t_{0}^{\mathrm{if}}=-\infty$, $t_{\tilde{K}+1}^{\mathrm{if}}=\infty$
and $\hat{K}$ the number of intersected interfaces. $\hat{K}$ is
unknown at the beginning and depends on the ray path.

\begin{rem}
Please note that determining the ray paths depends on the refractive
index $n$, which is the unknown quantity in our inverse problem.
Therefore, the actual ray paths are also unknown in the beginning
and are approached within the reconstruction process.\end{rem}
\begin{problem}
The inverse problem consists in determining the complex refractive
index 
\begin{equation}
\tilde{n}(x)=n(x)+i\alpha(x)\frac{c_{0}}{4\pi f},\qquad x\in\Omega,
\end{equation}
from the knowledge of the measured data $g_{\mathrm{abs}}(\varphi_{i},s_{j})=\ln\left(\frac{1}{\tau(\varphi_{i},s_{j})}\right)$
and $g_{\mathrm{ref}}(\varphi_{i},s_{j})=d(\varphi_{i},s_{j})$ of
a parallel scanning geometry with (inital) angles $\varphi_{i}=\frac{2\pi(i-1)}{p}$
for $i=1,...,p$ and (initial) distances $s_{j}=\frac{R}{q}j,$ for
$j=-q,...,q$. As we do not consider straight ray paths, we usually
have $g_{\mathrm{abs}}(-\varphi,s)\neq g_{\mathrm{abs}}(\varphi,s)$
or $g_{\mathrm{ref}}(-\varphi,s)\neq g_{\mathrm{ref}}(\varphi,s)$
in contrast to X-ray CT. So the angles $\varphi_{i}$ should be evenly
distributed over the interval $\left[0,\,2\pi\right)$. Our task is
to find the best possible approximation of $\tilde{n}$ satisfying
\begin{equation}
\begin{array}{ccc}
{\displaystyle \intop_{L(\varphi_{i},s_{j})}\left(n(x)-1\right)\,\mathrm{d}x=g_{\mathrm{ref}}(\varphi_{i},s_{j})} & \mathrm{and} & {\displaystyle \intop_{L(\varphi_{i},s_{j})}\alpha(x)\,\mathrm{d}x=g_{\mathrm{abs}}(\varphi_{i},s_{j})}\end{array}
\end{equation}
for $i=1,...,p$ and $j=-q,...,q$.

\newpage{}
\end{problem}

\section{Algorithm development\label{sec:algorithm-development}}

In the following, we develop a numerical algorithm for the calculation
of the complex refractive index $\tilde{n}(x)=n(x)+i\alpha(x)\frac{c_{0}}{4\pi f}$.
Since algorithms based on partial differential equations such as the
wave equation are very time-consuming, an efficient reconstruction
shall only consider dominating physical effects. Poor reconstructions,
which arise from the application of conventional reconstruction algorithms
from computed tomography to THz data, are often based on the missing
consideration of refraction and reflection losses. 

Therefore, we would like to modify ART for THz tomography. This is
done by using a priori information (i.~e.~the knowledge of interfaces)
and considering relevant physical laws (Snell's law and reflection
losses according to Fresnel equations). This results in a method, 
which takes refraction as well
as transmission into account and determines simultaneously the refractive index
and the absorption coefficient from measurements of path difference
and intensity data.

\subsection{ART}

Our method is based on the classical ART as described in
e. g. \cite[pp. 128, 137 ff.]{Natterer2001}. So we first give a brief
overview of ART. The main idea behind of this method is the application of Kaczmarz's
method to Radon's integral equation. We search for the solution of
the system
\begin{equation}
\intop_{L^{\nu}}f(x)\,\mathrm{d}x=g^{\nu},\qquad\nu=1,...,N,\qquad x\in\Omega\label{eq:ART-eq-integrals}
\end{equation}
where $L^{\nu}:=L(\varphi_{\nu},s_{\nu})$ are straight lines, $f$
a function (refractive index $f=n-1$ or absorption coefficient $f=\alpha$),
$g^{\nu}:=g(\varphi_{\nu},s_{\nu})$ the measured data (path difference
or intensity data) and $N=p(2q+1)$. Let 
\begin{equation}
\Omega\subset\bigcup_{\mu=1}^{M}S_{\mu}
\end{equation}
be a pixel discretization of the domain $\Omega$ with pixels $S_{\mu},\,\mu=1,...,M$.
We obtain a linear system by discretizing $f$ via 

\begin{equation}
f(x)\approx\tilde{f}(x)=\sum_{\mu=1}^{M}f_{\mu}\chi_{\mu}(x),\qquad x\in\Omega,\label{eq:f-discrete}
\end{equation}
with characteristic functions

\begin{equation}
\chi_{\mu}(x)=\left\{ \begin{array}{l}
1,\qquad\mathrm{if\,}x\in S_{\mu},\\
0,\qquad\mathrm{otherwise}.
\end{array}\right.
\end{equation}
The integral is thus approximated by 
\begin{equation}
\intop_{L^{\nu}}f(x)\,\mathrm{d}x\approx\sum_{\mu=1}^{M}a_{\nu\mu}f_{\mu}
\end{equation}
with $a_{\nu\mu}=\mathrm{length}\,(L^{\nu}\,\cap\,S_{\mu})$. Now,
we define $a^{\nu}:=\left(a_{\nu1},...,a_{\nu M}\right)^{T}\in\mathbb{R}^{M}$
and\linebreak{}
$F=(f_{1},f_{2},\text{...},f_{M})^{T}\in\mathbb{R}^{M}$. Using the
scalar product in $\mathbb{R}^{M}$, we obtain the discrete version
of (\ref{eq:ART-eq-integrals}), 
\begin{equation}
\left\langle a^{\nu},F\right\rangle =g^{\nu},\qquad\nu=1,...,N.\label{eq:ART-eq-discrete}
\end{equation}
We solve this linear system by applying Kaczmarz's method to (\ref{eq:ART-eq-discrete}).
The iteration then reads as

\begin{equation}
\begin{array}{lll}
F^{0} & := & \tilde{F}^{r}\\
F^{\nu} & = & {\displaystyle F^{\nu-1}+\lambda\frac{g^{\nu}-\left\langle a^{\nu},F^{\nu-1}\right\rangle }{\left\Vert a^{\nu}\right\Vert _{2}^{2}}a^{\nu},\qquad\nu=1,\,...,\,N}\\
\tilde{F}^{r+1} & := & F^{N}
\end{array}\label{eq:ART-iteration}
\end{equation}
with relaxation parameter $\lambda\in(0,2)$, iteration index $r=0,\,1,\,2,\,\text{...}$,
ray index $\nu$ and an arbitrary initial value $\tilde{F}^{0}$.

\subsection{Modified ART}

We modify ART taking refraction into account. This leads
to a more sophisticated calculation of the projection matrix $A=(a_{\nu\mu})_{\nu=1,...,N,\mu=1,...,M}$.
Furthermore, the consideration of Fresnel losses leads to a correction
of the measured intensity data $g_{\mathrm{abs}}$. Both linear systems
(for path difference $d$ and transmission coefficient $\tau$) show
a similar structure and can be solved simultaneously, using the same
ray paths for both models, by applying Kaczmarz's method to get the
complex refractive index.

\subsubsection{Consideration of Fresnel losses}

Reflection losses influence the measured transmission coefficient
$\tau$ by reducing the intensity of the transmitted signal in addition
to absorption losses. Neglecting reflection losses in the model will
thus lead to a systematic error in the calculation of the absorption
coefficient: the resulting values for $\alpha$ will be too high,
as a larger absorption coefficient corresponds to smaller transmission
coefficients. Therefore, our approach will be a correction of the
measured intensity data, so that only absorption losses are considered
for the reconstruction of the absorption coefficient. Mukherjee et~al.~\cite{Mukherjee2013}
introduced this approach successfully for cylindric objects. We expand
this approach to arbitrary interfaces and an initially unknown refractive
index.

The reflection loss $\rho$ of a single interface is given by equation
(\ref{eq:Fresnel}). For this purpose, we need the angle of incident
and refraction, which are known from the calculation of the ray paths
according to Snell's law, as well as the refractive indices on both
sides of the interface, given by (\ref{eq:n-for-both-sides-of-the-if}).
Let $\hat{K}\in\mathbb{N}$ be the number of interfaces intersecting
the ray path $L$ and $\rho_{k_{\iota}}$ be the reflection loss at
interface $\Xi_{k_{\iota}}$ with $k_{\iota}\in\left[1,\,K\right],\iota=1,...,\,\hat{K}$.
The emitted intensity is reduced due to the corresponding reflection
losses by the factor

\begin{equation}
C_{\mathrm{abs}}:=\prod_{\iota=1}^{\hat{K}}\left(1-\rho_{k_{\iota}}\right).
\end{equation}
In addition, the intensity of the ray is reduced by absorption losses
according to Lambert-Beer's law
\begin{equation}
a=\exp\left(-\int_{L}\alpha\,\mathrm{d}x\right).
\end{equation}
 As a consequence, the transmission coefficient 
\begin{equation}
\tau=a\prod_{\iota=1}^{\hat{K}}\left(1-\rho_{k_{\iota}}\right)
\end{equation}
is measured by the receiver. We are now able to correct the measured
transmission data (the right hand side of integral equation (\ref{eq:integral-eq-CT}))
by

\begin{equation}
\tilde{g}_{\mathrm{abs}}=\ln\left(\frac{I_{0}}{I}\right)=\ln\left(\frac{1}{a}\right)=\ln\left(\frac{1}{\tau}\prod_{\iota=1}^{\hat{K}}\left(1-\rho_{k_{\iota}}\right)\right)
\end{equation}
and obtain data that depends only on the absorption coefficient $\alpha$.

\subsubsection{Consideration of the refracted ray path}

A matrix entry $a_{\nu\mu}$ represents the length of the (refracted)
ray $\nu$ in pixel $\mu$ as in conventional ART. In contrast to
ART, we cannot calculate all entries $a_{\nu\mu}$ at once as the
path itself is influenced by the material parameters. Instead, we
have to determine first of all the parameter $t_{l}^{\mathrm{if}}$
for partial ray $L_{l}^{\nu}\left(\varphi_{l},s_{l}\right)$ (see
section \ref{sub:ray-path-depending-on-refraction}). Then, we calculate
$a_{\nu\mu}$ until reaching the nearest interface, determine the
new angle $\varphi_{l+1}$ according to Snell's law and the reflection
loss $\rho_{l}$ for the current interface $\Xi_{\tilde{k}}$ according
to Fresnel equations. Afterwards, these steps can be repeated for
the next partial ray $L_{l+1}^{\nu},\,l=1,...,\hat{K}$ and $\nu=1,...,N$.
A pseudo code for calculating the matrix $A$ and the reflection losses
is given in algorithm \ref{alg:calculating-matrix-A}.
\begin{algorithm}
\begin{algorithmic}[1]
\For{$\nu=1$ to $N$}
	\State{$\rho=0, \ t_{0}^{\mathrm{if}}=-\infty,\ t_{1}^{\mathrm{if}}=0,\  l=1$ }
	\While{$t_{l}^{\mathrm{if}}<\infty$}
		\State{Calculate  $t_{l}^{\mathrm{if}}$ by (\ref{eq:t-of-the-next-interface})}
		\State{Calculate $a_{\nu\mu}=\mathrm{length\,of}\,(L_{l}^{\nu}\cap S_{\mu})$ for $t\in\left[t_{l-1}^{\mathrm{if}},t_{l}^{\mathrm{if}}\right]$}
		\If{$t_{l}^{\mathrm{if}}<\infty$}		
			
			\State{Calculate partial ray $L_{l+1}^{\nu}(\varphi_{l+1},s_{l+1})$ by (\ref{eq:omega-l-plus-1-cases}) and (\ref{eq:s-l-plus-1-calculation}) }
			\State{Calculate reflection losses $\rho_{l}^{\nu}$ at $\Xi_{\tilde{k}}$ by (\ref{eq:Fresnel})}
			\State{$l=l+1$}
			 
		\EndIf
	\EndWhile
	\State{$\hat{K}^{\nu}=l-1;$} 
	\State{Calculate $C_{\mathrm{abs}}^{\nu}=\prod_{l=1}^{\hat{K}^{\nu}}\left(1-\rho_{l}^{\nu}\right)$}
\EndFor
\State{$\tilde{g}_{\mathrm{abs}}^{\nu}=\ln\left(\frac{C_{\mathrm{abs}}^{\nu}}{\tau^{\nu}} \right),\,\nu=1,...,N $}\\

\Return{$A, \tilde{g}_{\mathrm{abs}}$}

\end{algorithmic}

\caption{Calculating matrix $A$ and the corrected $\tilde{g}_{\mathrm{abs}}$\label{alg:calculating-matrix-A}}
\end{algorithm}

The reconstruction of $\mathrm{Re}(\tilde{n})=n$ has a direct influence
on the reconstruction of the absorption coefficient $\alpha$, because
the reflection losses are calculated by using the refractive indices
of the materials that share the interface. The remaining steps are
identical to ART. We can use the same matrix $A$ for the reconstruction
of both $n-1$ and $\alpha$ and thus the calculation of $F_{\mathrm{ref}}^{r}$
and $F_{\mathrm{abs}}^{r}$ takes place in the same loop, speeding
up the reconstruction significantly. In the worst case, the complex
refractive index is completely unknown and the initial guess is the
complex refractive index of air, meaning $F_{\mathrm{ref}}^{0}=F_{\mathrm{abs}}^{0}=0$.
In this case it makes sense to perform consecutively several applications
of Kaczmarz's method and use their approximations for an update of
the ray paths, calculating a new matrix $A$ and correction value
$C_{\mathrm{abs}}$ before the next iteration. This means, we calculate
the first approximation by conventional ART without considering the
interfaces and reflection losses and then the next steps follow with
modified ART based on this approximation. This update of ray paths
based on an initial guess like $F_{\mathrm{ref}}^{0}=0$ was discussed
before in \cite{Johnson1975} and \cite{Schomberg1978}, where it
was used in the context of reconstructions based on the eikonal equation.

\subsubsection{Numerical scheme \label{sub:numerical-scheme} }

\begin{figure}[H]
\begin{centering}
\includegraphics[width=0.85\textwidth]{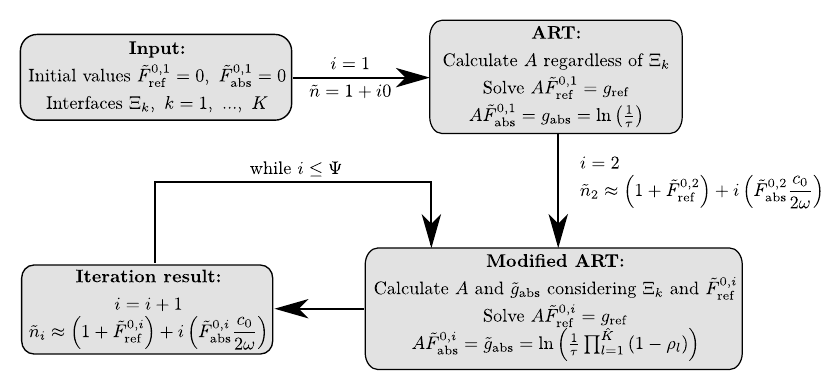}
\par\end{centering}

\caption{Numerical scheme of the modified ART\label{fig:numerical-scheme-modfied-ART}}
\end{figure}

Let $\Psi$ be the number of applications of Kaczmarz's method. We
can set the number of iterations $\psi=\left(r_{\mathrm{stop}}^{1},\text{...},r_{\mathrm{stop}}^{\Psi}\right)$
and the relaxation parameters for the refractive index $\lambda_{\mathrm{ref}}=\left(\lambda_{\mathrm{ref}}^{1},\text{...},\lambda_{\mathrm{ref}}^{\Psi}\right)$
and absorption coefficient $\lambda_{\mathrm{abs}}=\left(\lambda_{\mathrm{abs}}^{1},\text{...},\lambda_{\mathrm{abs}}^{\Psi}\right)$
for each iteration. An appropriate choice of $\psi$ is important,
due to the well-known phenomenon of the semiconvergence of iterative
methods for inverse problems. The stopping indices $r_{\mathrm{stop}}^{i},\,i=1,\,\text{...},\,\Psi,$
are the respective regularization parameters of the modified ART and
have to be chosen adequately (e.~g.~by the discrepancy principle,
see \cite{Defrise1987}). The algorithm for the calculation of an
approximation of the complex refractive index based on THz data is
shown in algorithm \ref{alg:modified-ART}. Altogether this leads
to the numerical scheme for the modified algebraic reconstruction
technique as outlined in figure \ref{fig:numerical-scheme-modfied-ART}.
\begin{algorithm}
\begin{algorithmic}[1]
\State{Inital values $\tilde{F}_{\mathrm{ref}}^{0,1},\ \tilde{F}_{\mathrm{abs}}^{0,1} $}
\For{$i=1$ to $\Psi$}
	\State{Calculate $A^{i}\left(\tilde{F}_{\mathrm{ref}}^{0,i}\right)$ and $\tilde{g}_{\mathrm{abs}}\left(\tilde{F}_{\mathrm{ref}}^{0,i}\right)$ by algorithm \ref{alg:calculating-matrix-A}}
	\State{$r=0$}
	\While{$r\leq r_{\mathrm{stop}}^{i}$}
		\State{$F_{\mathrm{ref}}^{0}=\tilde{F}_{\mathrm{ref}}^{r,i},\,F_{\mathrm{abs}}^{0}=\tilde{F}_{\mathrm{abs}}^{r,i}$}
		
			\For{$\nu=1$ to $N$}
			
				\State{ ${\displaystyle F_{\mathrm{ref}}^{\nu}=F_{\mathrm{ref}}^{\nu-1}+\lambda_{\mathrm{ref}}^{i}\frac{g_{\mathrm{ref}}^{\nu}-\left\langle a^{\nu},F_{\mathrm{ref}}^{\nu-1}\right\rangle }{\left\Vert a^{\nu}\right\Vert _{2}^{2}}a^{\nu}}$}
		
				\State{${\displaystyle F_{\mathrm{abs}}^{\nu}=F_{\mathrm{abs}}^{\nu-1}+\lambda_{\mathrm{abs}}^{i}\frac{\tilde{g}_{\mathrm{abs}}^{\nu}-\left\langle a^{\nu},F_{\mathrm{abs}}^{\nu-1}\right\rangle }{\left\Vert a^{\nu}\right\Vert _{2}^{2}}a^{\nu}}$}

			\EndFor
			
		\State{$\tilde{F}_{\mathrm{ref}}^{r+1}=F_{\mathrm{ref}}^{N},\,\tilde{F}_{\mathrm{abs}}^{r+1}=F_{\mathrm{abs}}^{N}$}
		\State{$r=r+1$}
				
	\EndWhile
	\State{$\tilde{F}_{\mathrm{ref}}^{0,i+1}=F_{\mathrm{ref}}^{N},\,\tilde{F}_{\mathrm{abs}}^{0,i+1}=F_{\mathrm{abs}}^{N}$} 
\EndFor\\
\Return{${\displaystyle \tilde{n}\approx\left(1+\tilde{F}_{\mathrm{ref}}^{0,\Psi+1}\right)+i\left(\tilde{F}_{\mathrm{abs}}^{0,\Psi+1}\frac{c_{0}}{2\omega}\right)}$}

\end{algorithmic}

\caption{Modified ART for calculating $\tilde{n}=n+i\kappa$\label{alg:modified-ART}}
\end{algorithm}
 
\section{Numerical experiments and results\label{sec:Numerical-experiments}}

In this section we would like to validate the modified ART using both
synthetic and real measured data. The results show an improvement
of the reconstruction obtained by the proposed technique in comparison
to conventional CT methods. The complex refractive index can be reconstructed
in more detail and with smaller errors by applying the new algorithm. 

First of all we take a look at reconstructions with synthetic and
noisy data. Afterwards, we compare the results obtained by modified
ART with those resulting from FBP and conventional ART on the basis
of real measured data. Synthetic and real measured data were both
recorded for $p=360$ initial angles and $141$ distances (meaning
$q=70$). Therefore one data set consists of $N=50760$ measurements.
Furthermore, the number of iterations $\psi$ and the number of applications
of  Kaczmarz's methods $\Psi$, i.~e.~the regularization parameters,
were chosen without a parameter-choice method for the sake of simplicity.
The iterations are just terminated when further iterations show no
significant improvement, but not yet divergence. 

Additionally, we use the a priori information on the position of the
interfaces to reset the refractive index outside the object to $n-1=0$
before each Kaczmarz's method.  Thus the angles of incidence and refraction
are calculated more accurately, as the ray is broken at an actual
discontinuity of the refractive index. The initial values $\ensuremath{\tilde{F}_{\mathrm{ref}}^{0,1}=\tilde{F}_{\mathrm{abs}}^{0,1}=0}$
are used for all reconstructions, therefore the first iteration is
equivalent to conventional ART.

\subsection{Synthetic data}

The first numerical example concerns a circle with a radius of $50\,\mathrm{mm}$,
a refractive index of $n=1.4$ and an absorption coefficient of $\alpha=0.05\mathrm{\,cm^{-1}}$
containing an embedded rectangle with dimensions of $25\,\mathrm{mm}\times20\,\mathrm{mm}$
and $n=1.7$, $\alpha=0.25\,\mathrm{cm^{-1}}$. Uniformly distributed
noise is added respectively to the synthetic transmission coefficient
and path difference data, such that the obtained noisy data $g^{\varepsilon}$
fulfills 
\begin{equation}
\left\Vert g-g^{\varepsilon}\right\Vert _{2}/\left\Vert g\right\Vert _{2}\leq0.05
\end{equation}
in both cases. 
\begin{figure}[h]
\hfill{}\subfloat[Refractive index $n(x),\,x\in\Omega$ ]{

\includegraphics[width=0.5\textwidth]{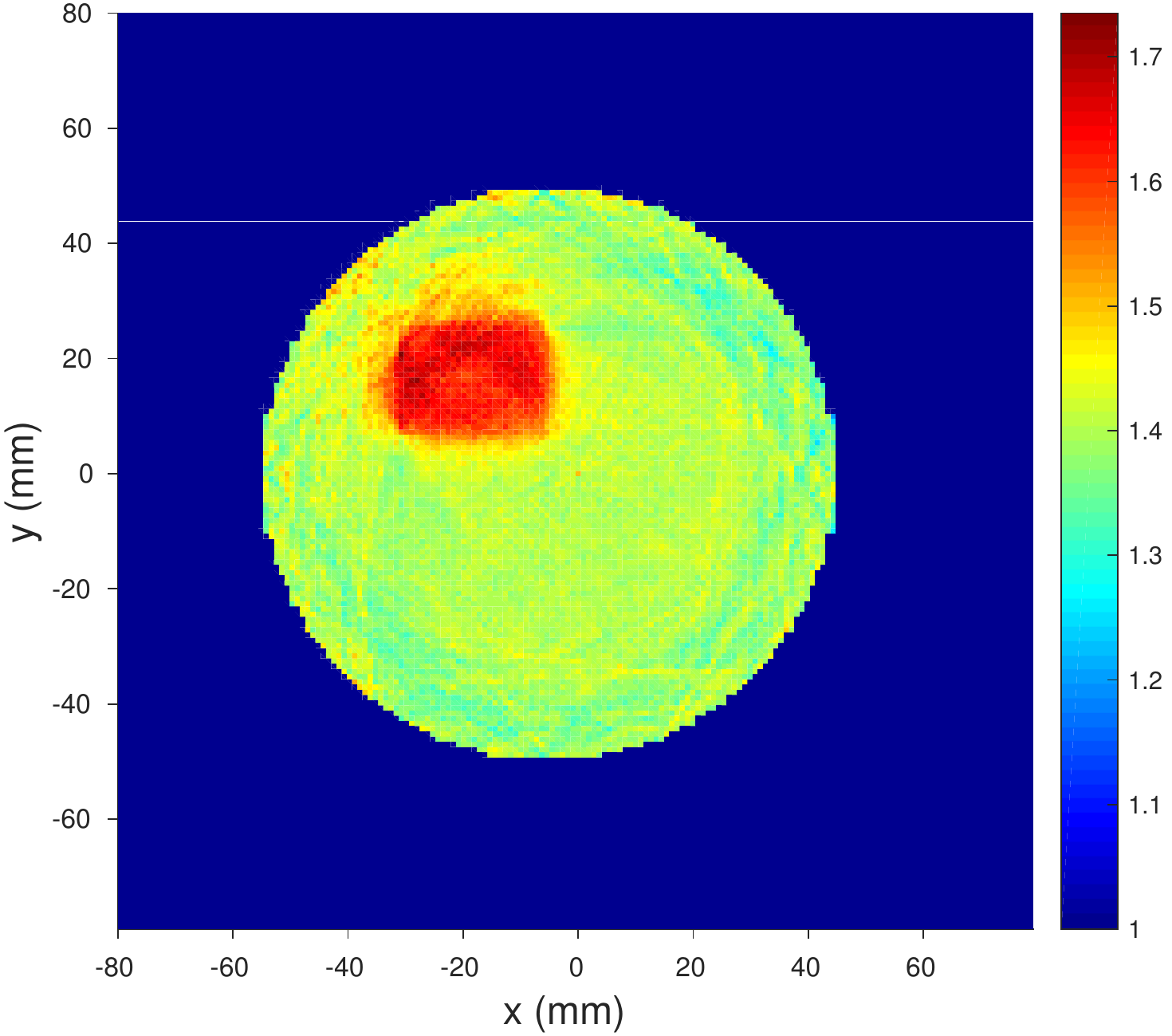}}\subfloat[Absolute error]{

\includegraphics[width=0.5\textwidth]{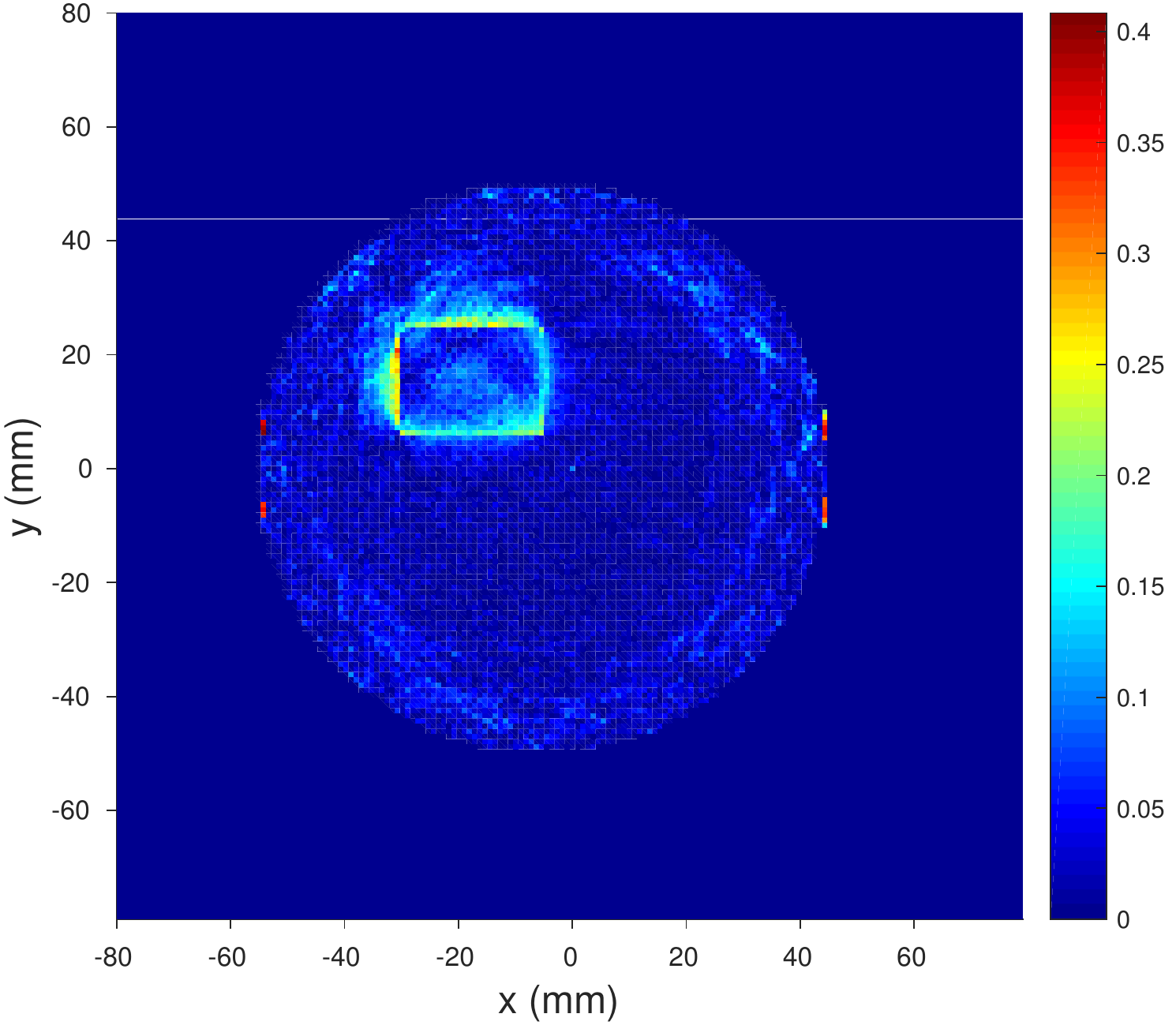}}\hfill{}\caption{Reconstruction of the refractive index based on synthetic data\label{fig:rec-circle-ref}}
\end{figure}
As there is no information about $\tilde{n}$, we perform several
applications of Kaczmarz's method for the reconstruction as described
in section \ref{sub:numerical-scheme}. Here we used $\Psi=5$ iterations
with the following parameters:

\begin{equation}
\begin{array}{c}
\psi=\left(3,\,3,\,5,\,7,\,5\right)\\
\lambda_{\mathrm{ref}}=\left(0.01,\,0.01,\,0.006,\,0.002,\,0\right)\\
\lambda_{\mathrm{abs}}=\left(0.002,\,0.004,\,0.004,\,0.004,\,0.003\right)
\end{array}
\end{equation}
We perform only a reconstruction of the absorption coefficient $\alpha$
in last iteration ($\lambda_{\mathrm{ref}}^{5}=0$). This approach
is reasonable in so far as we use the approximation of the previous
iteration step of the refractive index $n$ for calculating the reflection
losses for the reconstruction of $\alpha$. Therefore, we can expect
an improvement of $\alpha$ in an additional iteration step. 
\begin{figure}[h]
\hfill{}\subfloat[Absorption coefficient $\alpha(x),\,x\in\Omega$ ]{

\includegraphics[width=0.5\textwidth]{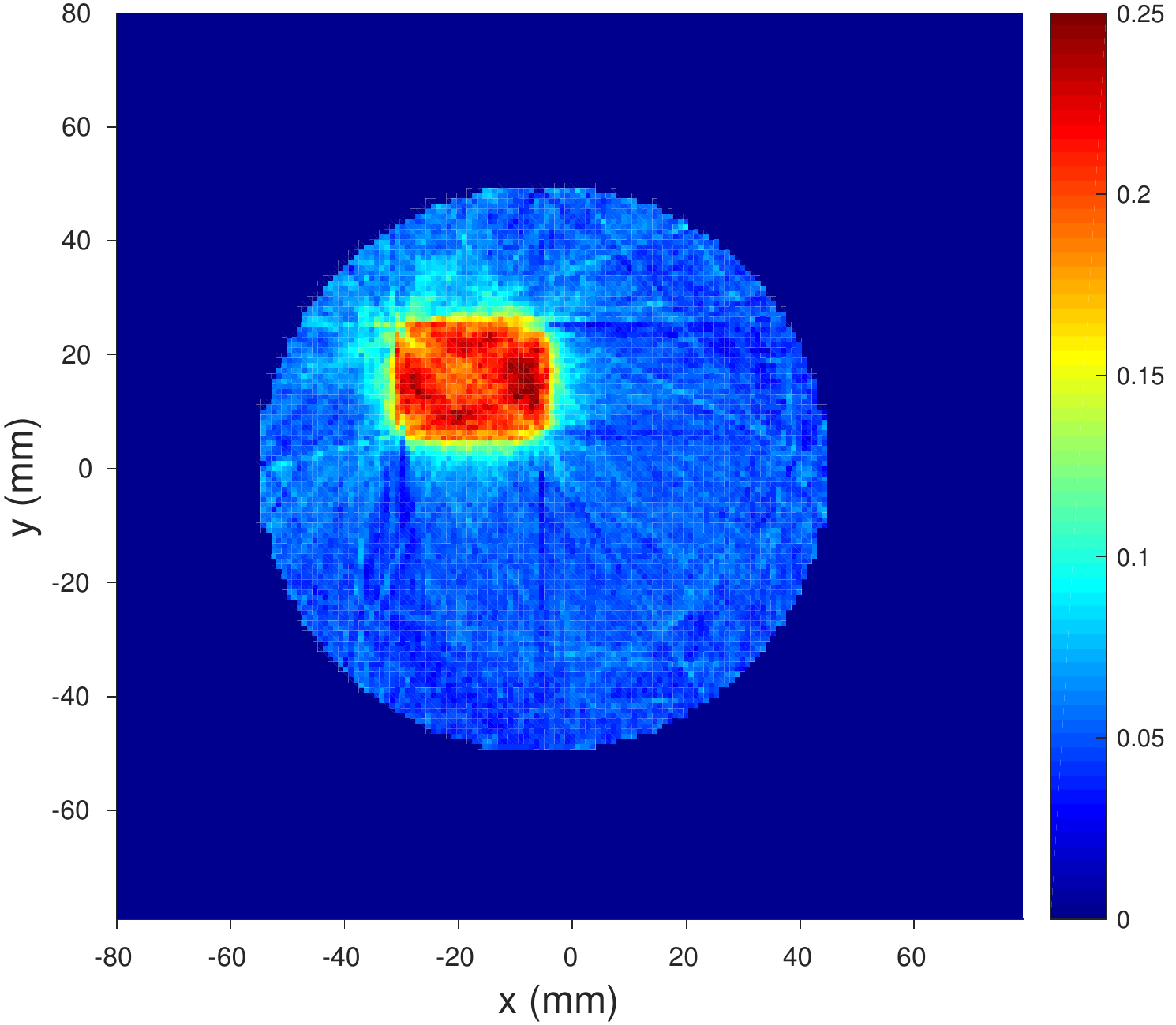}}\subfloat[Absolute error]{

\includegraphics[width=0.5\textwidth]{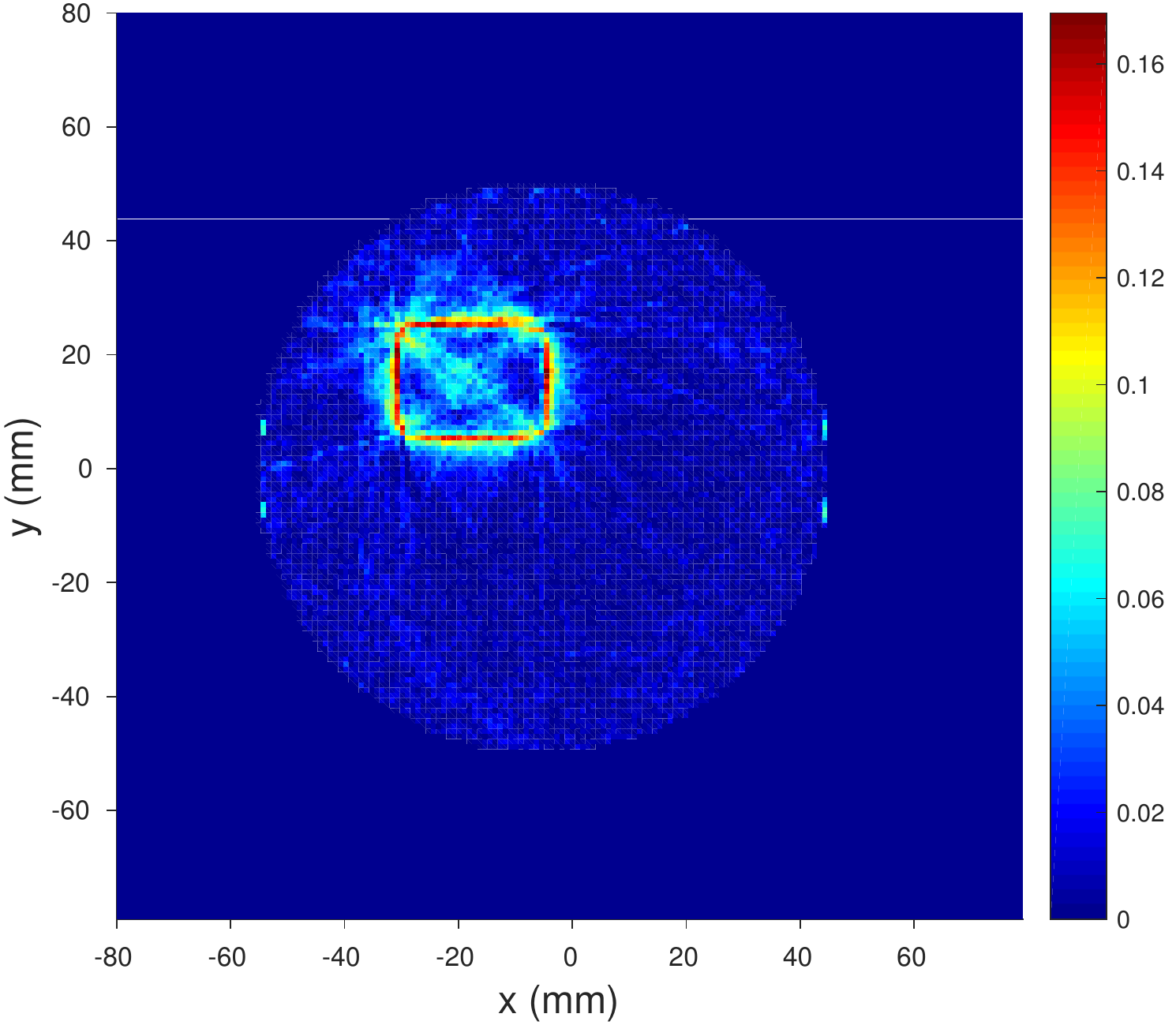}}\hfill{}\caption{Reconstruction of the absorption coefficient based on synthetic data\label{fig:rec-circle-abs}}
\end{figure}
Figures \ref{fig:rec-circle-ref} and \ref{fig:rec-circle-abs} show
the reconstructions of the refractive index and the absorption coefficient
with the modified ART. The results are good in spite of the
large jump discontinuities of the refractive index. The
absolute error is small in large parts of the object (see figures
on the right), whereas the the biggest local errors appear near the
interfaces and corners. Hereby, the prerequisite is ensured that modified
ART provides correct results with synthetic data. As we now have to
deal with ray paths that can be arbitrary refracted at interfaces,
it is possible that some pixels are less frequently intersected by
the ray paths than others, which causes additional artifacts.

\subsection{Comparison with conventional reconstruction methods using real data}

The next step is to compare modified ART with conventional ART and
FBP based on real measured data. In this context, we first give a
brief description of the used THz tomograph and then move on to the
reconstructions.

\subsubsection{THz tomography system}

The measured data were acquired at \emph{SKZ}\footnote{~The German Plastics Center in W\"urzburg, \href{http://www.skz.de/}{http://www.skz.de/}}.
A THz tomography system was designed, which works on the basis of
the FMCW (Frequency Modulated Continuous Wave) principle. Here the
frequency of the THz signal is a sawtooth wave in the range of $70$
to $110\,\mathrm{GHz}$ and allows measurements of intensity and phase.
The latter results in the delayed travel time or path difference data.
Dispersion can be neglected in this operating frequency range and
the complex refractive index can be assumed to be constant \cite{Lamb1996,Piesiewicz2007}.
Thus the transmission coefficient and the path difference are averaged
over this frequency band. The optical system consists of two lenses
for each receiver and transmitter and is similar to the architecture
of a Kepler telescope. The Gaussian beam is focused on the object
with a beam waist of about $8\,\mathrm{mm}$. 
\begin{rem}
The generation of synthetic data always provides suitable values $g_{\mathrm{abs}}$
and $g_{\mathrm{ref}}$, regardless of how strong the rays are refracted.
However, the THz system uses just one emitter and one detector, and
in the case of strong refraction the transmitted ray may miss the
detector, resulting in a transmission coefficient of $\tau=0$. Better
results can often be achieved by sorting out measured data of rays
with $\tau\leq\varepsilon_{\mathrm{miss}}$ and $\varepsilon_{\mathrm{miss}}$
sufficiently small, meaning neglecting them completely for the iterations.
It must be taken into account that this case also occurs for strong
absorbers and then leads to a distortion of the reconstruction results. 
\end{rem}

\subsubsection{Block with multiple refractive indices}

\begin{figure}[h]
\centering{}\includegraphics[width=0.45\textwidth]{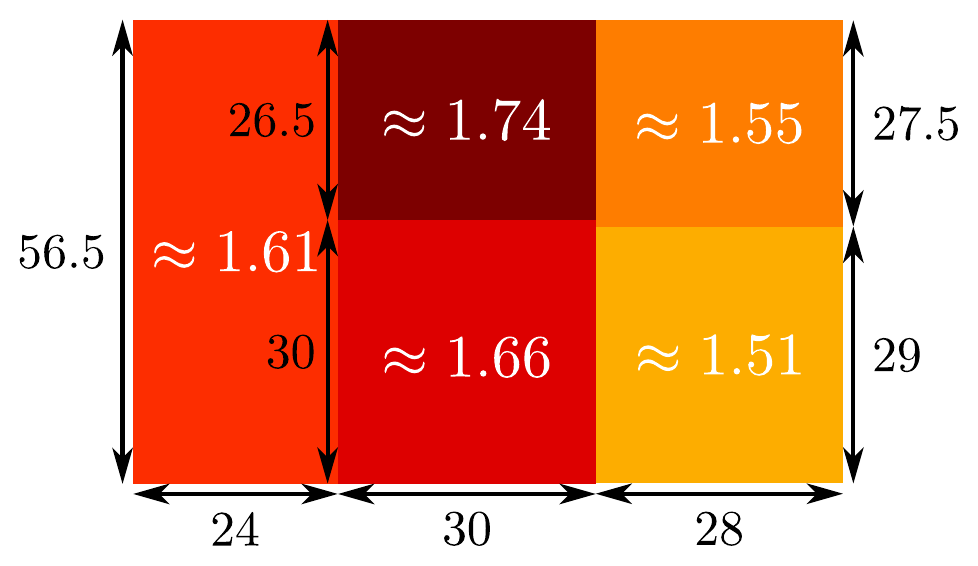}\caption{Refractive indices and proportions (in $\mathrm{mm}$) of the glued-together
block \label{fig:ref-ind-and-prop-of-block}}
\end{figure}

\begin{figure}[h]
\hfill{}

\subfloat[Refractive index $n(x),\,x\in\Omega$ ]{

\includegraphics[width=0.5\textwidth]{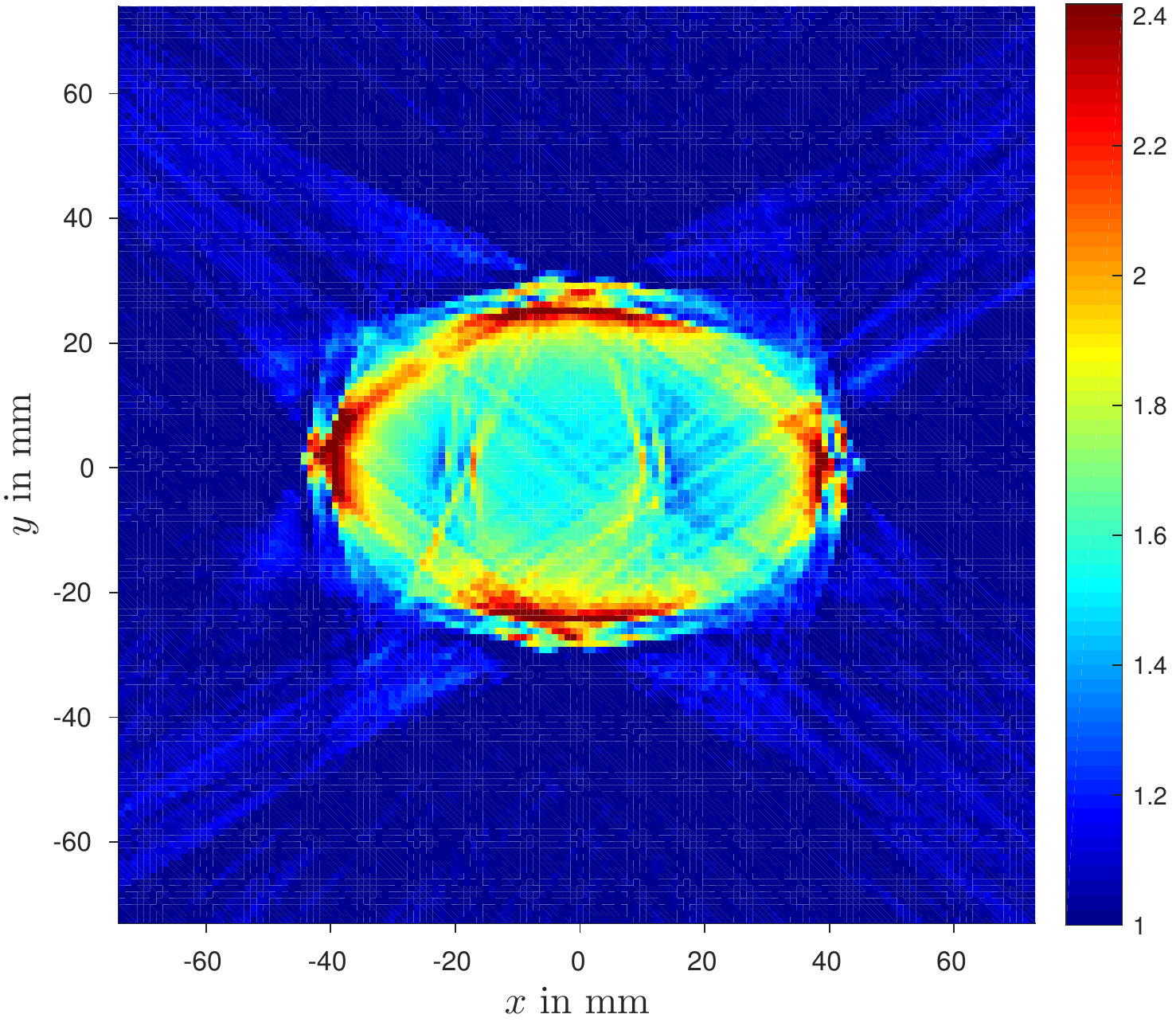}}\subfloat[Absolute error]{

\includegraphics[width=0.5\textwidth]{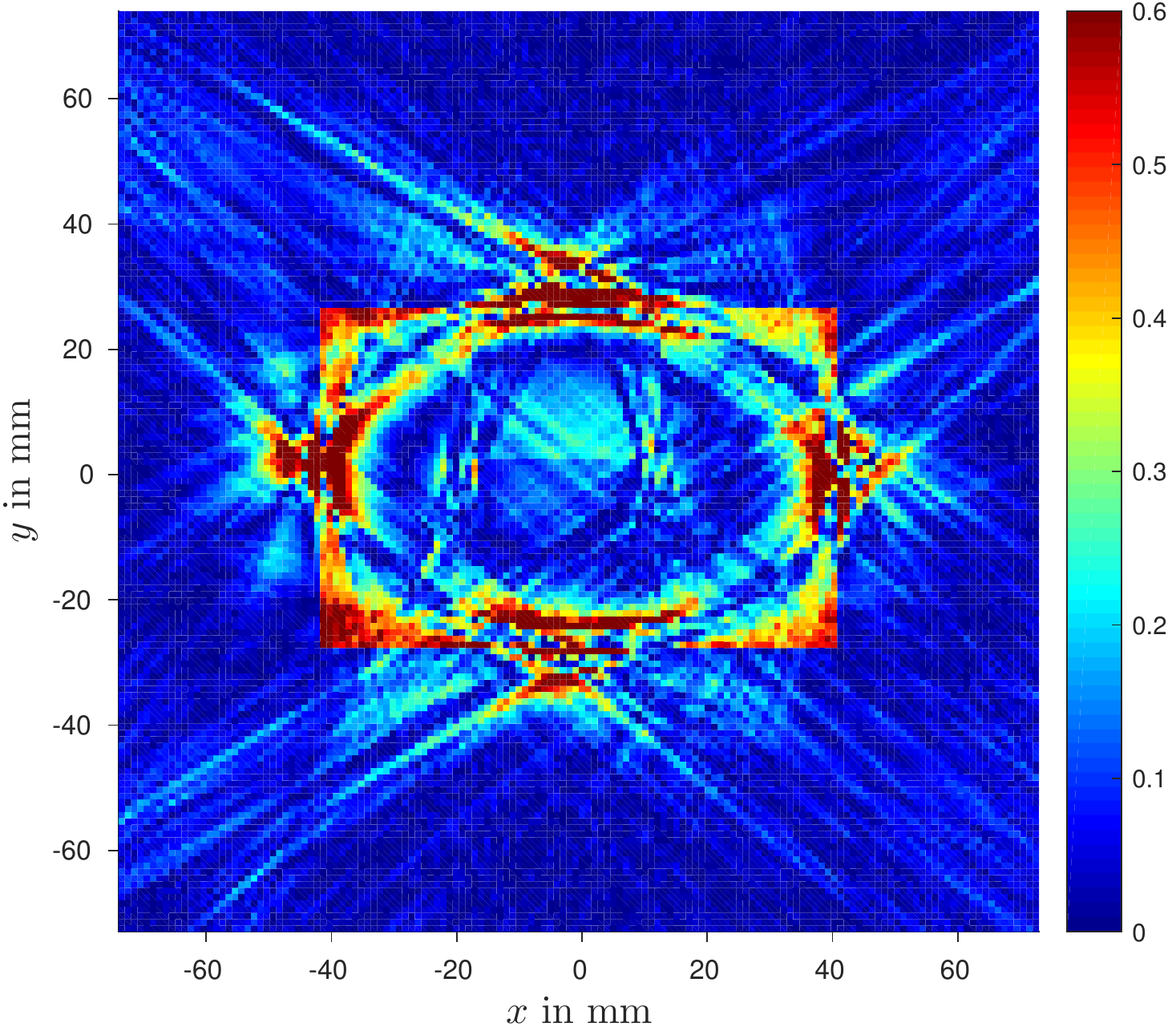}}\hfill{}\caption{Filtered backprojection: reconstruction of the refractive index with
real measured data \label{fig:FBP-rec-5-block}}
\end{figure}

\begin{figure}[H]
\hfill{}\subfloat[Refractive index $n(x),\,x\in\Omega$ ]{

\includegraphics[width=0.5\textwidth]{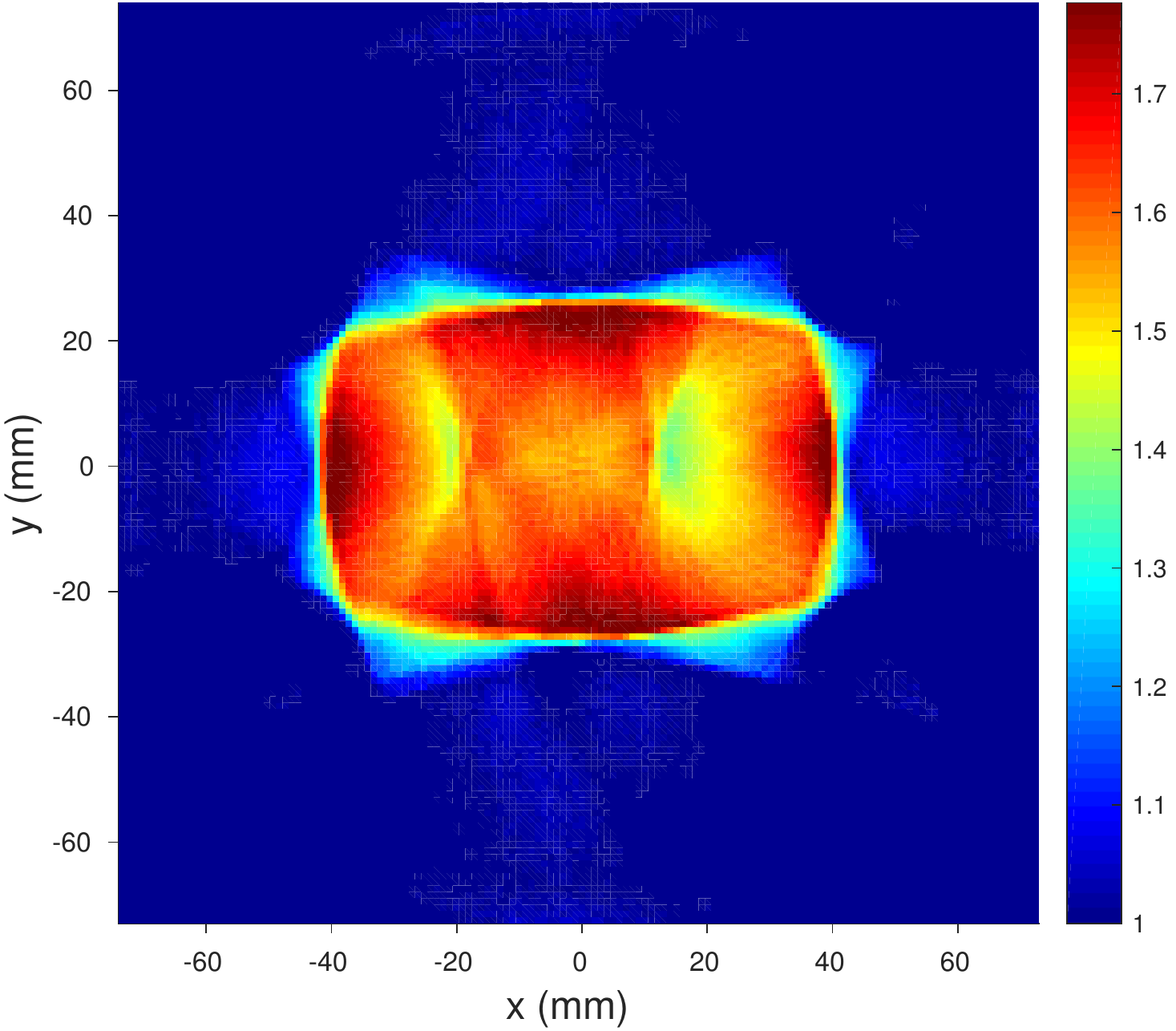}}\subfloat[Absolute error]{

\includegraphics[width=0.5\textwidth]{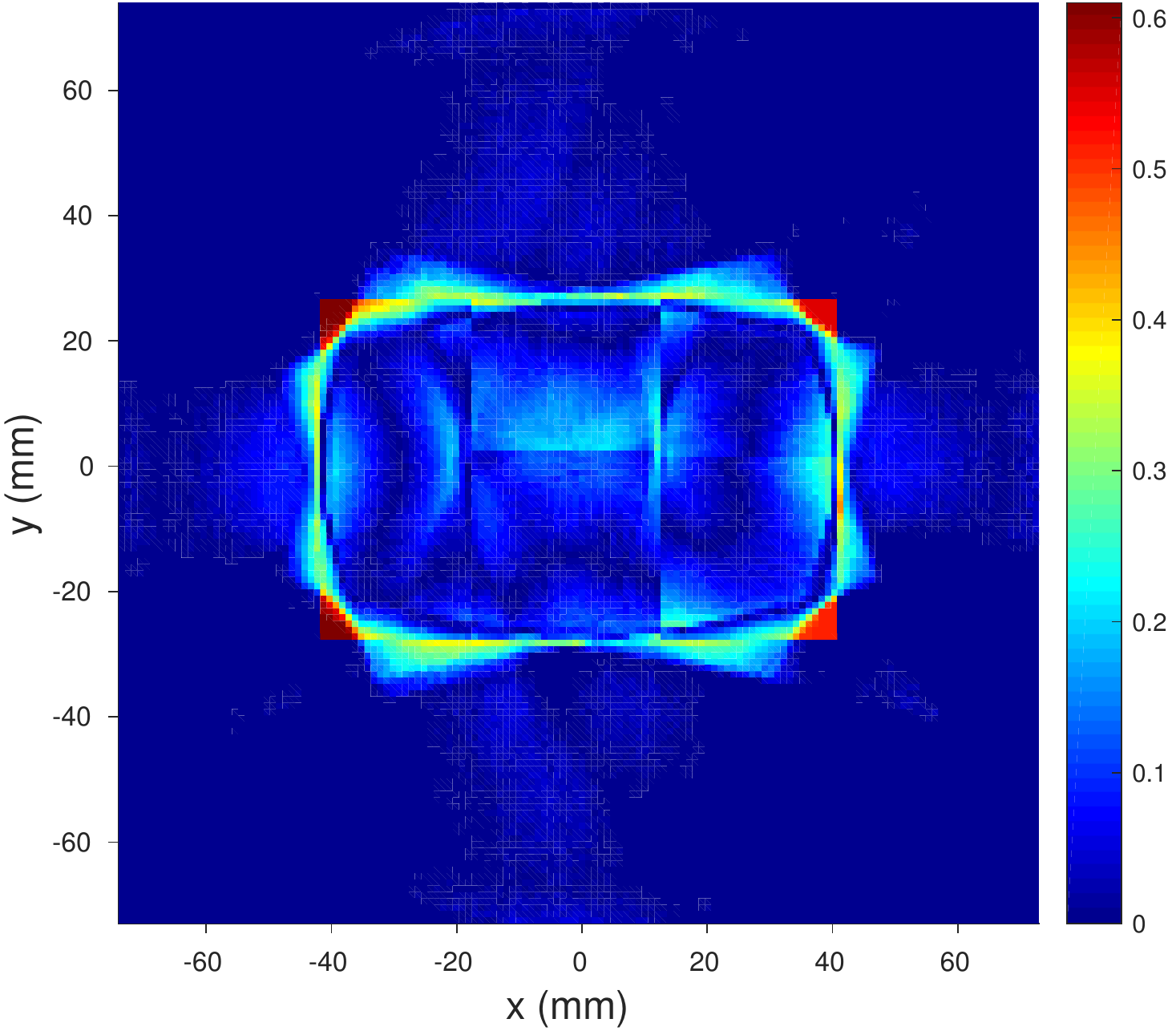}}\hfill{}\caption{ART: reconstruction of the refractive index with real measured data
\label{fig:ART-rec-5-block}}
\end{figure}

\begin{figure}[H]
\hfill{}\subfloat[Refractive index $n(x),\,x\in\Omega$ ]{

\includegraphics[width=0.5\textwidth]{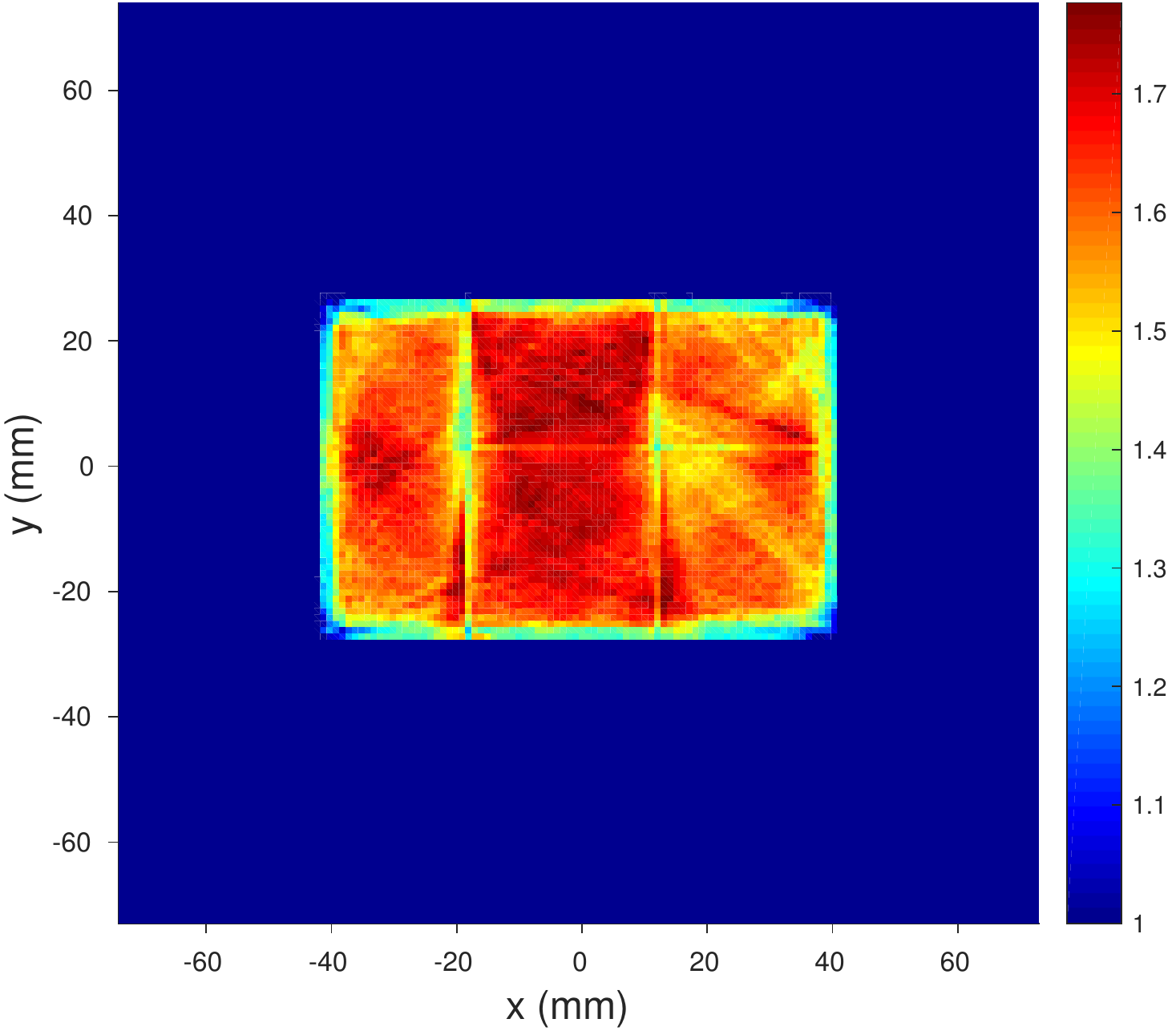}}\subfloat[Absolute error]{

\includegraphics[width=0.5\textwidth]{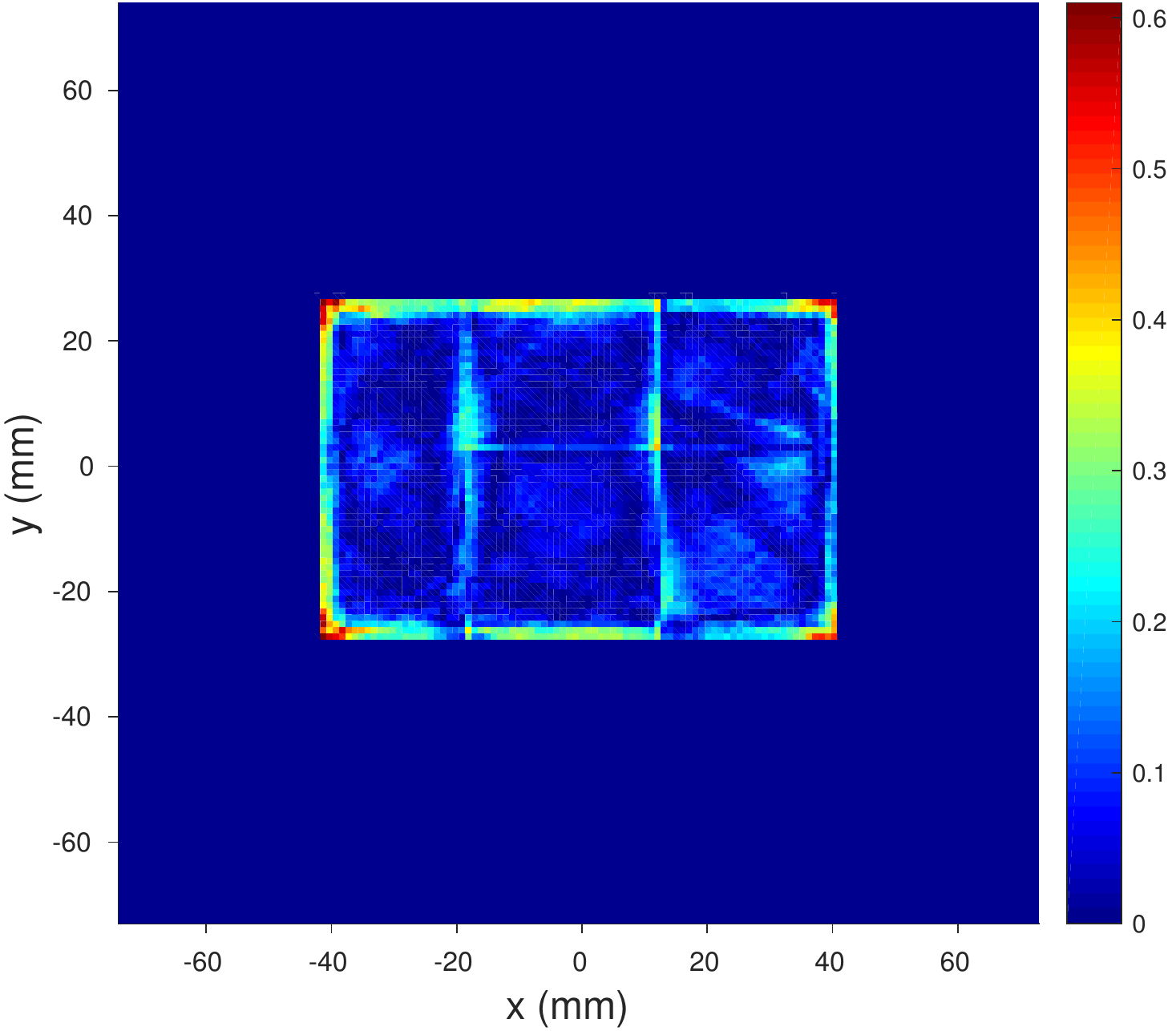}}\hfill{}\caption{Modified ART: reconstruction of the refractive index with real measured
data\label{fig:Modified-ART-rec-5-block}}
\end{figure}

The next step is to validate our modified ART with a more complex
object and real measured data. Furthermore, we would like to compare
the filtered backprojection, conventional ART and modified ART. For
this purpose a test object consisting of different plastic blocks
with different refractive indices that where glued together forming
a block was scanned. Figure \ref{fig:ref-ind-and-prop-of-block} shows
the measures and the refractive indices of the object.
The conventional ART no longer achieves significant improvements after
15 iterations with $\lambda_{\mathrm{ref}}=0.005$. The results of
modified ART are obtained with $\psi=\left(3,\,3,\,5,\,12\right)$
and $\lambda_{\mathrm{ref}}=\left(0.005,\,0.01,\,0.0075,\,0.003\right)$.
THz rays resulting in transmission coefficients with $\tau\leq0.05$
are ignored in both algebraic methods. Furthermore, we used the filtered
backprojection with the Shepp-Logan filter and $\varepsilon_{\mathrm{miss}}=0$.
The application of FBP to THz data leads to a very bad reconstruction
(see figure \ref{fig:FBP-rec-5-block}) with very high values of $n$
at the outer interfaces. The shape of the object is completely lost
and cannot be reconstructed.
 The reconstruction of the conventional ART is slightly better (see
figure \ref{fig:ART-rec-5-block}). But the outer shape of the object
also appears too curved and the reconstruction of the refractive index
produces very high values especially around the interfaces.
 Conversely, the reconstruction with the modified ART is much better,
see figure \ref{fig:Modified-ART-rec-5-block}, and good results can
also be achieved for more complex objects. The modified ART outperforms
the conventional methods FBP and ART, especially around corners, regarding
shape and size of the refractive index. This is also evident from
the comparison of the three methods regarding the absolute error (see
respective figures on the right-hand side). 
 The absolute error of the modified ART is almost everywhere smaller.
The biggest deviations appear again at the corners and interfaces.
Consequently, our modified ART can reconstruct more complex objects
with smaller errors and in more detail, while the individual blocks
are hardly detectable with ART. The differentiation of the various
plastics is difficult with modified ART as well, because the values
of the refractive indices are close to each other. But we can clearly
distinguish between polyamide (middle top, with the highest refractive
index) and polyethylene (right bottom, with the lowest refractive
index). Finally, the reconstructions with modified ART are clearly
superior to those of FBP and ART, showing that refraction and reflection
losses play an important role in THz tomography.

\section{Conclusion and future work}

In this paper, we have presented a new, efficient reconstruction algorithm
for THz tomography and have validated it with the help of synthetic
and real measured data. Our modified ART takes refraction and reflection
losses into account and reconstructs the complex refractive index
based on path difference and transmission coefficient measurements.
Therefore we use a refracted ray path according to Snell's law. Furthermore,
modified ART is a hybrid algorithm exploiting two different tomographic
principles (travel time and transmission tomography) and the respective
data simultaneously.

Numerical experiments show that the reconstruction of material parameters
and interfaces provide more accurate results than conventional CT
algorithms like ART. We have already pointed out in the introduction
that finding a balance between an efficient reconstruction and a detailed
modeling is important. Hence, the modified ART considers only the
main physical effects like refraction, absorption and reflection losses.
The advantage of this approach is that we can avoid the expensive
solution of partial differential equations like wave equation or Helmholtz
equation. As we neglect further physical effects in our modeling,
our proposed method has the following limitations: Since we do not
consider diffraction effects, objects containing small, non-absorbing
defects that lead to strong diffraction are poorly reconstructed.
Although the first results are promising, the reconstruction algorithm
for THz tomography can be improved by the following approach: An important
property of THz radiation is the propagation according to Gaussian
optics. At this point, we neglect this physical property, suggesting
the consideration of the Gaussian beam profile as a promising next
step. Taking into account the Gaussian beam profile, we might also
consider cases in which only a part of the transmitted beam hits the
receiver and thus improve the reconstruction of the absorption coefficient.

To reduce the time required for reconstruction by parallelization,
using the simultaneous algebraic reconstruction technique (SART, \cite{Andersen1984})
or more general block algebraic iterative methods (see e. g. \cite{Soerensen2014})
is conceivable. Since the modifications to conventional ART concern
only the matrix $A$ and the data $g_{\mathrm{abs}}$, the iteration
of e. g. SART is directly applicable, too.

\section*{Acknowledgments}
This project was funded by The German Federation of Industrial Research Associations
(AiF) under 457 ZN which is gratefully acknowledged by the authors.

\bibliographystyle{siam}
\bibliography{gIPE_THz_ART}

\end{document}